\documentclass[11pt,twoside]{article}

\usepackage[latin1]{inputenc}

\usepackage[english]{babel}
\usepackage{amsfonts,amssymb,amsmath}
\usepackage{mathrsfs}
\usepackage{mathtext}
\usepackage{textcomp}
\usepackage{graphicx}
\usepackage{epsfig}
\usepackage{rotating}
\usepackage{tikz}
\usepackage{fancybox}

\usetikzlibrary{arrows,decorations.pathmorphing,backgrounds,fit}

\usepackage[toc,page]{appendix}

\usepackage[dvips, colorlinks=true,pdfstartview=FitV, linkcolor=blue,citecolor=red,urlcolor=blue]{hyperref}

\newcommand{\lra}{\longrightarrow}
\newcommand{\setm}{\setminus}
\newcommand{\ra}{\rightarrow}
\newcommand{\DS}{\displaystyle}

\newcommand{\C}{\mathbb{C}}
\newcommand{\R}{\mathbb{R}}
\newcommand{\Z}{\mathbb{Z}}

\newcommand{\N}{\mathbb{N}}

\renewcommand{\S}{\mathbb{S}}
\renewcommand{\P}{\mathrm{P}}
\newcommand{\T}{\mathrm{T}}
\newcommand{\Tt}{\mathcal{T}}
\newcommand{\A}{\mathcal{A}}
\newcommand{\Ah}{\hat{\mathcal{A}}}
\newcommand{\Her}{\mathbf{H}}

\newcommand{\J}{\mathbb{J}}

\newcommand{\Sig}{\Sigma}
\newcommand{\di}{\mathbf{d}}
\newcommand{\Id}{\mathrm{Id}}
\newcommand{\inter}{\mathrm{int}}
\newcommand{\leng}{\mathbf{leng}}
\newcommand{\orth}{\mathbf{orth}}
\newcommand{\pQ}{(\mathcal{Q})}

\newcommand{\card}{\mathrm{Card}}

\renewcommand{\Re}{\mathrm{Re}}
\renewcommand{\Im}{\mathrm{Im}}
\newcommand{\Tr}{\mathrm{Tr}}

\newcommand{\ie}{\textit{i.e. }}
\renewcommand{\geq}{\geqslant}
\renewcommand{\leq}{\leqslant}
\newcommand{\vide}{\varnothing}

\newcommand{\Hg}{\mathcal{H}(k_1,\dots,k_n)}
\newcommand{\Hgi}{\mathcal{H}_1(k_1,\dots,k_n)}

\newcommand{\transpose}[1]{{\vphantom{#1}}^{\mathit t}{#1}}

\newtheorem{RecThm}{Theorem}

\newcommand{\HomgA}{\mathrm{Homeo}^+(S_g,\hat{\mathcal{A}})}
\newcommand{\HomgAo}{\mathrm{Homeo}^+_0(S_g,\hat{\mathcal{A}})}

\newcommand{\GgA}{\Gamma(S_g,\hat{\mathcal{A}})}

\newcommand{\Tgn}{\mathcal{T}(g,n)}
\newcommand{\Ggn}{\Gamma(g,n)}

\newcommand{\Homgno}{\mathrm{Homeo}^+_0(S_g,\{p_1,\dots,p_n\})}

\newcommand{\Tsn}{\mathcal{T}(0;n)}

\newcommand{\Gsn}{\Gamma(0;n)}

\newcommand{\Teich}{\mathcal{T}_\T(\bar{\alpha};\bar{\beta})}

\newcommand{\TETt}{\widetilde{\mathcal{T}}^{\mathrm{et}}(\Ah,\underline{\alpha})^*}
\newcommand{\TET}{\mathcal{T}^{\mathrm{et}}(\Ah,\underline{\alpha})}
\newcommand{\MET}{\mathcal{M}^{\mathrm{et}}(\Ah,\underline{\alpha})}

\newcommand{\TETv}{\mathcal{T}^{\mathrm{et}}(\Ah,\underline{\alpha})^*}
\newcommand{\TETvi}{\mathcal{T}_1^{\mathrm{et}}(\Ah,\underline{\alpha})^*}
\newcommand{\METv}{\mathcal{M}^{\mathrm{et}}(\Ah,\underline{\alpha})^*}
\newcommand{\METvi}{\mathcal{M}_1^{\mathrm{et}}(\Ah,\underline{\alpha})^*}

\newcommand{\TSSv}{\mathcal{T}(\S^2,\underline{\alpha})^*}
\newcommand{\TSS}{\mathcal{T}(\S^2,\underline{\alpha})}

\newcommand{\TSSvi}{\mathcal{T}_1(\S^2,\underline{\alpha})^*}
\newcommand{\MSS}{\mathcal{M}(\S^2,\underline{\alpha})}
\newcommand{\MSSi}{\mathcal{M}_1(\S^2,\underline{\alpha})}
\newcommand{\MSSv}{\mathcal{M}(\S^2,\underline{\alpha})^*}
\newcommand{\MSSvi}{\mathcal{M}_1(\S^2,\underline{\alpha})^*}
\newcommand{\TRSn}{\mathcal{TR}(\S^2,\{p_1,\dots,p_n\})}

\newcommand{\STt}{\mathbf{S}_\mathcal{T}}

\newcommand{\U}{\mathcal{U}}
\newcommand{\V}{\mathcal{V}}
\newcommand{\VTt}{\mathrm{V}_\mathcal{T}}

\newcommand{\ATt}{\mathbf{A}_\mathcal{T}}
\newcommand{\ATtA}{\mathbf{A}_{\mathcal{T},\mathcal{A}}}

\newtheorem{theorem}{Theorem}[section]
\newtheorem{definition}[theorem]{Definition}
\newtheorem{proposition}[theorem]{Proposition}
\newtheorem{prop-def}[theorem]{Proposition-definition}
\newtheorem{corollary}[theorem]{Corollary}
\newtheorem{lemma}[theorem]{Lemma}

\newcommand{\dem}{\noindent {\textbf{Proof:} }}
\newcommand{\rem}{ \noindent \textbf{Remark: }}

\newcommand{\carre}{\hfill $\Box$\\}

\begin{document}

\title{\huge  TRIANGULATIONS AND VOLUME FORM ON MODULI SPACES OF FLAT SURFACES} 

\author{ {\large DUC-MANH NGUYEN} \\ Laboratoire de Mathématiques, Bât.425\\ Université Paris Sud XI\\ 91405 Orsay Cedex, France.}

\date{}

\maketitle

\begin{abstract}
In this paper, we study the moduli spaces of flat surfaces with cone singularities verifying the following property: there exists a union of
disjoint geodesic tree on the surface such that the complement is a translation surface. Those spaces can be viewed as deformations of the
moduli spaces of translation surfaces in the space of flat surfaces. We prove that such spaces are quotients of flat complex affine manifolds by
a group acting properly discontinuously, and preserving a parallel volume form. Translation surfaces can be considered as a special case of flat
surfaces with erasing forest, in this case, it turns out that our volume form coincides with the usual volume form (which are defined via the
period mapping) up to a multiplicative constant. We also prove  similar results for the moduli space of flat metric structures on the
$n$-punctured sphere with prescribed cone angles up to homothety. When all the angles are smaller than $2\pi$, it is known (cf. \cite{Thu}) that
this moduli space is a complex hyperbolic orbifold. In this particular case, we prove that our volume form induces a volume form which is equal
to the complex hyperbolic volume form up to a multiplicative constant.

\end{abstract}


\tableofcontents

\section{Introduction}

While the moduli space of (half)-translation surfaces has been studied for a rather long time, little is known about the topology of its strata,
we only know that each stratum contains finitely many (in fact, at most three) connected components thank to recent works of Kontsevich-Zorich
\cite{KonZo}, and Lanneau \cite{Lann}. To investigate the topology of the strata, it is often useful to study their deformations. In this paper,
we study some deformations of the moduli space of translation surfaces in the space of flat surfaces with conical singularities. For this
purpose, we deform a translation surface by replacing a single singular point by a cluster of singular points with fixed cone angles, the angles
are chosen so that this replacement can be carried out in a neighborhood  of the former singular point without changing the geometric structure
on its complement. For technical reasons, we also assume that the new singular points are the vertices of  a geodesic tree which will be called
an {\em erasing tree}, the union of such trees is called an {\em erasing forest}.\\

\noindent On a flat surface with erasing forest, vector fields which are invariant by parallel transport can be defined on the complement of the
erasing forest, we will call such fields {\em parallel vector fields}.\\

\noindent As all the trees in an erasing forest shrink to points, a flat surface with erasing forest becomes a usual translation surface,
therefore, the moduli space of surface with erasing forest can be viewed as a deformation of the moduli space of translation surfaces. Note
that, by definition, a translation surface is a particular flat surface with erasing forest where all the erasing trees are points. \\

\noindent On a flat surface of genus zero, a geodesic tree whose vertex set is the set of  singular points is automatically an erasing tree.
Since such a tree always exists, the study of flat surface with erasing forest provides us with a common treatment for translation surfaces and
flat surface of genus zero. Before getting into formal definitions and precise statements, let us resume the main results of this paper in the
following

\begin{RecThm}

The moduli space of flat surfaces with erasing forest and unitary parallel vector field (here, the erasing forest consists of trees which are
isomorphic to those of some fixed family of topological trees, the cone angles at the vertices of the erasing forest are fixed) is a quotient of
a flat complex affine manifold by a group acting properly discontinuously, preserving a parallel volume form $\mu_{\mathrm{Tr}}$.\\

\noindent In the case of translation surfaces, up to a multiplicative constant, the volume form $\mu_\Tr$ agrees with the usual one, which is
defined by the period mapping. \\

\noindent For the case of flat surfaces of genus zero, the volume form $\mu_\Tr$ induces a volume form $\hat{\mu}_\Tr^1$ on the moduli space of
flat surfaces of genus zero of unit area with fixed cone angles. If all the cone angles are less than $2\pi$, then $\hat{\mu}^1_\Tr$ agrees with
the complex hyperbolic volume form defined by Thurston up to a multiplicative constant.\\

\end{RecThm}

In a forthcoming paper, we will prove that the volume of the moduli space of flat surfaces with erasing forest with respect to $\mu_\Tr$,
normalized by some energy functions involving the area of the surface and the length of the trees in the erasing forest, is finite. Using this
result, one can recover the classical results of Masur-Veech, and Thurston on the finiteness of the volume of the moduli space of translation
surface, and of the moduli space of flat surface of genus zero.\\

\noindent \textbf{Acknowledgements:} The author is deeply indebted to François Labourie for suggesting the subject, and for the guidance without
which this work would not be accomplished. He would like also to warmly thank Pascal Hubert and the referee for numerous useful comments on this
work.

\section{ Background and statement of main results}

\subsection{Flat surfaces and their moduli spaces}

Let $\Sig$ be an oriented, connected, closed surface, and $Z$ a finite subset of $\Sig$. A {\em flat surface structure on $\Sig$ with conical
singularities at $Z$} is an Euclidean metric structure on $\Sig \setm Z$ so that for every $s$ in $Z$, a neighborhood of $s$ is isometric to an
Euclidean cone. Note that if the cone angle at a point in $Z$ is $2\pi$ then this point is actually  a {\em regular} point. Throughout this
paper, we will call a flat surface with conical singularities a {\em flat surface}.\\

\noindent {\em Translation surfaces} are flat surfaces such that the holonomy of any closed curve, which does not contain any singular point, is
a translation of $\R^2$. On such surfaces, given a direction $\theta$ in $\S^1$, there exists a foliation $\mathcal{F}_\theta$ of the surface by
parallel lines in direction $\theta$.\\

The moduli spaces of flat surfaces have been studied by numerous authors \cite{BavGh}, \cite{EskMas}, \cite{Thu}, \cite{Vee93}. For the case of
polyhedral flat surfaces, let $\alpha_1,\dots,\alpha_n$ be $n$ positive real numbers such that $\DS{\alpha_1+\dots+\alpha_n=2\pi(n-2)}$. Let
$C(\alpha_1,\dots,\alpha_n)$ denote the space of all flat surface structures of unit area on the sphere $\S^2$ having exactly $n$ singular
points with cone angles $\alpha_1,\dots,\alpha_n$. The following result is proved in \cite{Thu}.

\begin{RecThm}[Thurston]
Assume that all the angles $\alpha_i, \; i=1,\dots,n$ are smaller than $2\pi$, then $C(\alpha_1,\dots,\alpha_n)$ is a complex hyperbolic
orbifold of dimension $n-3$, whose volume is finite.

\end{RecThm}

A necessary condition for a flat surface to be a translation surface is that the cone angle at every singular point must belong to $2\pi\N$. As
a consequence, there are no translation surface structures on the $2$-sphere.\\

A translation surface together with a foliation by parallel lines can be identified to a pair $(M,\omega)$, where $M$ is a closed connected
Riemann surface, and $\omega$ is a holomorphic $1$-form on $M$. By this identification, a zero of order $k$ of the holomorphic $1$-form
corresponds to a singular point with cone angle $2\pi(k+1)$.\\

\noindent Given positive integers $g,k_1,\dots,k_n$ such that $\DS{k_1+\dots+k_n=2g-2}$, let $\Hg$ denote the set of all pairs $(M,\omega)$,
where $M$ is a Riemann surface of genus $g$, and $\omega$ is a holomorphic $1$-form on $M$ which has exactly $n$ zeros with orders
$k_1,\dots,k_n$. We denote by $\Hgi$ the subset of $\Hg$ which corresponds to the set of translation surfaces of area one. It is well known (cf.
\cite{EskOk}, \cite{Kon}, \cite{Vee90}, \cite{Zor}) that $\Hg$ is a complex algebraic orbifold of dimension $2g+n-1$, and there exists a volume
form $\mu_0$ on $\Hg$ which is induced by the {\em period mapping} $\Phi$. The map $\Phi$ is defined locally by the integrals of $\omega$ over a
family of curves representing a basis of $H_1(M,Z(\omega), \Z)$, where $Z(\omega)$ is the set of zeros of $\omega$. This map sends a
neighborhood of $(M,\omega)$ in $\Hg$ to an open subset of $\C^{2g+n-1}$. The volume form $\mu_0$ is the pull-back under $\Phi$ of the Lebesgue
measure on $\C^{2g+n-1}$.\\

\subsection{Erasing forest}

Given a flat surface $\Sig$, a {\em tree} in $\Sig$ is the image of an embedding from a topological tree into $\Sig$. We consider an isolated
point as a special tree which has only one vertex and no edges. A {\em forest} in $\Sig$ is a union of disjoint trees in $\Sig$. A tree in
$\Sig$ is said to be {\em geodesic } if each of its edges is a geodesic segment in $\Sig$. A forest is said to be {\em geodesic} if it is a
union of geodesic trees.\\

\begin{definition}[Erasing forest]
\label{defET}

Let $\Sig$ be a compact connected flat surface without boundary. An {\em erasing forest}  on $\Sig$ is a  geodesic forest whose vertex set
contains all the singular points of $\Sig$ such that, if $c$ is a closed curve in $\Sig$ which does not intersect any tree in the forest, then
the holonomy of $c$ is a translation of $\R^2$.

\end{definition}

\rem \begin{itemize}

\item[$\bullet$] If an erasing forest contains only one tree, we will call this tree an {\em erasing tree}.

\item[$\bullet$] Since points are trees, a translation surface has an obvious erasing forest which is the union of all singular points.

\item[$\bullet$] If $\Sig$ is a flat surface homeomorphic to the $2$-sphere, and $A$ is a geodesic tree on $\Sig$ which connects all the
singular points, then $A$ is automatically an erasing tree, since $\Sig\setm A$ is a topological disk. Throughout this paper, we will call a
flat surface of genus zero, a {\em spherical flat surface}. On any flat surface, there always exists a geodesic trees whose vertex set is
exactly the set of singular points (cf. Proposition \ref{ETprEx}), therefore, spherical flat surfaces can be considered, not in a unique way, as
special flat surfaces with an erasing tree.

\end{itemize}

Given a flat surface $\Sig$ with an erasing forest $\hat{A}$, we call a non vanishing vector filed defined on the complement of $\hat{A}$  a
{\em parallel vector field} if its integral lines are parallel lines in the local charts of the Euclidean metric structure.\\

\subsection{Main results}

We fix two integers $g\geq 0$, $n>0$, such that $2g+n-2>0$, and positive real numbers $\alpha_1,\dots,\alpha_n$ verifying
$\alpha_1+\dots+\alpha_n=2\pi(2g+n-2)$. In what follows, $S_g$  will be a fixed compact, oriented, connected flat surface of genus $g$, without
boundary. We also assume that there exists a geodesic erasing forest $\Ah=\sqcup_{i=1}^m\A_i$ on $S_g$, where each $\A_i$ is a geodesic tree.
Let $p_1,\dots,p_n$ denote the vertices of the trees in $\Ah$, and assume that the cone angle at $p_i$ is $\alpha_i$. Let $\V$ denote the set
$\{p_1,\dots,p_n\}$. Recall that, by definition, all the singular points of $S_g$ are contained in $\V$, but some of the points $p_i$ may be
regular.\\

\noindent Let $\HomgA$ denote the group of orientation preserving homeomorphisms of $S_g$ which fix the points in $\V$, and preserve the forest
$\Ah$. Let $\HomgAo$ be the normal subgroup of $\HomgA$ consisting of all elements which can be connected to  $\Id_{S_g}$ by an isotopy fixing
the points in $\V$.

\begin{definition}[Mapping class group preserving a forest] \label{ETdefMCG}

The {\em mapping class group} of the pair $(S_g,\Ah)$ is the quotient group

$$\GgA=\HomgA/\HomgAo.$$

\end{definition}

We start by defining the space of flat metric structures having an erasing forest isomorphic to $\Ah$ together with a marking, \ie the
Teichmüller space of flat surfaces with erasing forest, we then identify the moduli space of flat surfaces with erasing forest with the quotient
of this Teichmüller space under the action of the group $\GgA$.\\

\noindent Let $\underline{\alpha}$ denote the set $\{\alpha_1,\dots,\alpha_n\}$, and let $\TETt$ denote the collection of all pairs
$(\Sig,\phi)$, where $\Sig$ is an oriented flat surface of genus $g$, and $\DS{\phi: S_g \ra \Sig}$ is an orientation preserving homeomorphism
verifying

\begin{itemize}
\item[(i)] $\phi$ maps the set $\{p_1,\dots,p_n\}$ bijectively to the set of singularities of $\Sig$ such that the cone angle at the point
$\phi(p_i)$ is $\alpha_i, \; i=1,\dots,n$.

\item [(ii)] The image of $\Ah$ under $\phi$ is an erasing forest on $\Sig$.

\end{itemize}

\noindent We define an equivalence relation on $\TETt$ as follows: two pairs $(\Sigma_1,\phi_1)$ and $(\Sigma_2,\phi_2)$ are equivalent if there
exists an isometry $h: \Sig_1\longrightarrow \Sig_2$ such that the homeomorphism $\phi_2^{-1}\circ h\circ\phi_1$ is an element of $\HomgAo$. The
equivalence class of a pair $(\Sig,\phi)$ is denoted by $[(\Sig,\phi)]$. Let $\TETv$ be the set of equivalence classes of this relation.\\

\noindent Clearly, the group $\GgA$ acts on $\TETv$, the quotient space $\TETv/\GgA$ will be called the {\em moduli space of flat surfaces with
marked erasing forest} and denoted by $\METv$. We denote  by $\TETvi$ the set of equivalence classes $[(\Sig,\phi)]$, where $\Sig$ is a flat
surface of area one, and by $\METvi$ the quotient space $\TETvi/\GgA$.

\begin{definition}
\label{ETdefTeich}

The {\em Teichmüller space of flat surfaces with marked erasing forest and parallel vector field} is the set of equivalence classes of triples
$(\Sig,\phi,\xi)$, where $(\Sig,\phi)$ is a pair in $\TETt$, and $\xi$ is a unitary parallel vector field on $\Sig\setm\phi(\Ah)$, and the
equivalence relation is defined as follows: $(\Sig_1,\phi_1,\xi_1)$ and $(\Sig_2,\phi_2,\xi_2)$ are equivalent if there exists an isometry
$h:\Sig_1 \ra \Sig_2$ such that

\begin{itemize}
\item[$\bullet$] $h_*\xi_1=\xi_2$,

\item[$\bullet$] $\phi_2^{-1}\circ h \circ \phi_1 \in \HomgAo$.

\end{itemize}

We denote this space by $\TET$, the equivalence class of a triple $(\Sig,\phi,\xi)$ will be denoted by $[(\Sig,\phi,\xi)]$.\\

\noindent The {\em moduli space of flat surfaces with marked erasing forest and parallel vector field} is the quotient space $\TET/\GgA$, we
denote it by $\MET$.

\end{definition}

A point in $\MET$ is then a triple $(\Sig,\hat{A},\xi)$, where $\Sig$ is a flat surface having exactly $n$ singularities, with cone angles
$\alpha_1,\dots,\alpha_n$, $\hat{A}$ is an erasing forest on $\Sig$ isomorphic to $\Ah$, and $\xi$ is a unitary parallel vector field on
$\Sig\setm\hat{A}$. Two topological trees are isomorphic if there exists a continuous mapping from one to the other which induces bijections on
the two sets of vertices, and the two set of edges. Note that, by definition, we have a numbering on the set of vertices of $\hat{A}$, and
hence, a numbering on the set of edges of $\hat{A}$.\\

Recall that a {\em flat complex affine manifold} is a manifold with an atlas whose transition maps are complex affine transformations. We can
now state the main results of this paper.

\begin{theorem}\label{ETthB}
The space $\TET$ is a flat complex affine manifold of dimension

\begin{itemize}
\item[$\bullet$] $2g+n-1$ if $\alpha_i\in 2\pi\N$ for every $i\in\{1,\dots,n\}$.

\item[$\bullet$] $2g+n-2$ otherwise.
\end{itemize}

\noindent The group $\GgA$ acts properly discontinuously on $\TET$, and preserves a parallel volume form which will be denoted by $\mu_\Tr$.\\

\end{theorem}

In the case where $\Ah$ is a union of points, which implies that $S_g$ is a translation surface, the space $\MET$ can be identified to a stratum
$\Hg$ for some appropriate integers $k_1,\dots,k_n$. Recall that, on $\Hg$, we have a volume form $\mu_0$ which is defined by  the period
mapping. We have

\begin{proposition}\label{ETprC}
 On each connected component of $\Hg$, there exits a constant $\lambda$ such that $\mu_\mathrm{Tr}=\lambda \mu_0$.
\end{proposition}

Let us now focus on the case $g=0$, in this case, we identify $S_g$ to the standard sphere $\S^2$, and fix $n$ points $p_1,\dots,p_n$ on $\S^2$
with $n\geq 3$.  Fix  $n$ positive real numbers $\underline{\alpha}=(\alpha_1,\dots,\alpha_n)$ such that $\alpha_1+\dots+\alpha_n=2\pi(n-2)$.
The {\em Teichmüller space of spherical flat surfaces} having $n$ singularities with cone angles $\alpha_1,\dots,\alpha_n$ is the set of
equivalence classes of pairs $(\Sig,\phi)$, where

\begin{itemize}
\item[$\bullet$] $\Sig$ is a spherical flat surface having $n$ singularities with cone angles $\alpha_1,\dots,\alpha_n$.

\item[$\bullet$] $\phi$ is a homeomorphism from $\S^2$ to $\Sig$, which sends $\{p_1,\dots,p_n\}$ onto the set of singularities of $\Sig$ such
that the cone angle at $\phi(p_i)$ is $\alpha_i$.

\item[$\bullet$] The equivalence class of $(\Sig,\phi)$ is the set of all pairs $(\Sig,\phi')$, where $\phi'$ is a homeomorphism isotopic to
$\phi$ by an isotopy which is constant on the set $\{p_1,\dots,p_n\}$.
\end{itemize}

\noindent We denote this Teichmüller space by $\TSSv$. The equivalence class of a pair $(\Sig,\phi)$ in $\TSSv$ is denoted by $[(\Sig,\phi)]$.
Let $\TSSvi$ denote the subset of $\TSSv$ consisting of surfaces of unit area.\\

Let $\Gsn$ denote the modular group of the punctured sphere $\S^2\setm\{p_1,\dots,p_n\}$. Let $\MSSv$, and $\MSSvi$ denote the quotients
$\TSSv/\Gsn$, and $\TSSvi/\Gsn$ respectively.\\

Let $\TSS$ denote the product space $\TSSv\times\S^1$. We extend the action of $\Gsn$ onto $\TSS$ by assuming that  $\Gsn$ acts trivially on the
$\S^1$ factor of $\TSS$. Let $\MSS$, and $\MSSi$ denote the spaces $\MSSv\times \S^1$, and $\MSSvi\times \S^1$ respectively. First, we have

\begin{proposition}\label{SSprA}
$\TSS$ can be endowed with a flat complex affine manifold structure of dimension $n-2$, on which $\Gsn$ acts properly discontinuously.\\
\end{proposition}

By Proposition \ref{ETprEx}, we know that, on any flat surface of genus zero, there always exists an erasing tree whose vertex set is the set of
singularities. Therefore, we can identify a neighborhood of a point $([(\Sig,\phi)],e^{\imath\theta})$ in $\TSS$ to a neighborhood of a point
$[(\Sig,\phi,\xi)]$ in $\TET$, where $\Ah$ contains only one tree. By this identification, we get a volume form on a neighborhood of
$([(\Sig,\phi)],e^{\imath\theta})$, which is induced by the volume form $\mu_\mathrm{Tr}$, depending on the choice of the erasing tree. Our main
result in this case is the following

\begin{theorem}\label{SSthB}
Let $A_1$ and $A_2$ be two erasing trees on $\Sig$, and let $\mu_{A_1},\mu_{A_2}$ denote the volume forms corresponding to $A_1,A_2$
respectively which are defined on a neighborhood of $([(\Sig,\phi)],e^{\imath\theta})$ in $\TSS$, then we have

$$\mu_{A_1}=\mu_{A_2}.$$

\noindent Consequently, we get a well-defined volume form $\mu_\mathrm{Tr}$ on $\TSS$ which is invariant under the action of $\Gsn$.

\end{theorem}

The volume form $\mu_\mathrm{Tr}$ induces naturally a volume form $\hat{\mu}^1_\mathrm{Tr}$ on $\MSSvi$. In the case where $0<\alpha_i<2\pi$,
for $i=1,\dots,n$, according to the Thurston's result, $\MSSvi$ can be equipped with a complex hyperbolic metric structure, hence, we have a
volume form $\mu_{\mathrm{Hyp}}$ induced by this metric on $\MSSvi$. In this situation, we have

\begin{proposition}\label{SSprC}

There exists a constant $\lambda$ such that $\hat{\mu}^1_\mathrm{Tr}=\lambda\mu_\mathrm{Hyp}$.

\end{proposition}

This paper is organized as follows: from Section \ref{ETAdTrSec} to Section \ref{ETVolSec} we prove Theorem \ref{ETthB}, the proof of
Proposition \ref{ETprC} is given in Section \ref{prfCSec}, and the proof of Proposition \ref{SSprA} in Section \ref{SSAffStSec}, Section
\ref{SSVolSec} is devoted to the proof of Theorem \ref{SSthB}, and finally the proof of Proposition \ref{SSprC} is given in Section
\ref{SSCompSec}.\\

\section{Admissible triangulation}\label{ETAdTrSec}

Let $[(\Sig,\phi,\xi)]$ be a point in  $\TET$. Following the method of Thurston in \cite{Thu}, we will construct local charts of $\Teich$ about
$[(\Sig,\phi,\xi)]$ using geodesic triangulations of $\Sig$. In preparation for this construction, we first define

\begin{definition}[Admissible triangulation]
\label{defAT} Let $(\Sig,\phi)$ be a pair in $\TETt$, an {\em admissible triangulation} of $(\Sig,\phi)$ is a triangulation $\T$ of $\Sig$ such
that:

\begin{itemize}
\item[$\bullet$] The set of vertices of $\T$ is the set $V=\phi(\{p_1,\dots,p_n\})$.

\item[$\bullet$] Every edge of $\T$ is a geodesic segment.

\item[$\bullet$] The erasing forest $\phi(\Ah)$ is contained in the 1-skeleton of $\T$.

\end{itemize}

\end{definition}

The aim of this section is to show the existence and uniqueness up to isotopy of such triangulations on $\Sig$. First, we have

\begin{proposition}[Existence of admissible triangulations] \label{ETprB1}
There always exists an admissible triangulation $\T$ for any pair $(\Sig,\phi)$ in $\TETt$.
\end{proposition}

\dem If we cut the surface $\Sig$ along the edges in the forest $\phi(\Ah)$, then we will obtain a flat surface with piece-wise geodesic
boundary. Therefore, the proposition is a particular case of the following lemma

\begin{lemma}
Let $\hat{\Sig}$ be a flat surface with piece-wise geodesic boundary if the boundary is not empty, and $\hat{V}$ be a finite subset of
$\hat{\Sig}$ which contains all the singular points. Then there exists a geodesic triangulation of $\hat{\Sig}$ whose vertex set is $\hat{V}$.
\end{lemma}

\noindent This is a well-known fact, for a proof, the reader can see for example \cite{BobSpr}, \cite{Vee93}, or \cite{EskMas}. \carre

Next, we have

\begin{proposition}[Uniqueness of admissible triangulations up to isotopy]\label{ETprB2}
Let $\T_1$ and $\T_2$ be two admissible triangulations of $(\Sig,\phi)$. If there exists a homeomorphism $h$ in $\HomgAo$ such that
$h(\T_1)=\T_2$ then $\T_1=\T_2$.

\end{proposition}

\dem This proposition is a direct consequence of Lemma \ref{ETlmA} below.

\begin{lemma}\label{ETlmA}

Let $\Sig$ be a flat surface without boundary. Let $V=\{x_1,\dots,x_n\}$ be a finite subset of $\Sig$ such that $\Sig\setm V$ contains only
regular points, and suppose that $\chi(\Sig\setm V)<0$. Let $\gamma$ and $\gamma'$ be two simple geodesic arcs of $\Sig$ having the same
endpoints in $V$ (the two endpoints may coincide). Assume that $\gamma$ and $\gamma'$ are homotopic with fixed endpoints relative to $V$, then
we have $\gamma\equiv\gamma'$.

\end{lemma}

\dem We first observe that there exist no Euclidean structures on a closed disk such that its boundary is the union of two geodesic segments.
This is just a consequence of the Gauss-Bonnet Theorem.\\

\noindent Since $\chi(\Sig\setm V) <0$, the universal covering of $\Sig\setm V$ is the open disk $\Delta=\{z \in \C :\; |z| < 1\}$. The flat
metric structure on $\Sig\setm V$ give rise to a flat metric structure on $\Delta$ (which is  not complete). Now, let $\tilde{\gamma}$ be a lift
of $\gamma$  in $\Delta$ whose endpoints are contained in the boundary of $\Delta$ . By lifting the homotopy from $\gamma$ to $\gamma'$, we get
a lift $\tilde{\gamma}'$ of $\gamma'$ which has the same endpoints as $\tilde{\gamma}$. Note that by assumption, $\tilde{\gamma}$ and
$\tilde{\gamma}'$ are two geodesic in $\Delta$.\\

\noindent The two curves $\tilde{\gamma}$ and $\tilde{\gamma}'$ may have intersections, but in any case, we can find (at least) an open disk $D$
which is bounded by two arcs, one is a subsegment of $\tilde{\gamma}$, the other is a subsegment of $\tilde{\gamma}'$. Consequently, the open
disk $D$ is isometric to the interior of an Euclidian disk which is bounded by two geodesic segments. Since such a disk cannot exist, we have
$\tilde{\gamma}\equiv \tilde{\gamma}'$, and the lemma follows.\carre

\section{Flat complex affine structure}\label{ETAffStSec}

In this section, we prove that the space $\TET$ is a flat complex affine manifold and compute its dimension. Let $\mathcal{TR}(S_g,\Ah)$ denote
the set of all equivalence classes of triangulations (not necessarily geodesic) of $S_g$ which contain the forest $\Ah$, and whose vertex set is
$\V$, two triangulations are equivalent if they are isotopic relative to $\V$. Let $\Tt$ be an element of $\mathcal{TR}(S_g,\Ah)$, we consider
$\Tt$ as a particular triangulation of $S_g$, and denote by $\U_{\Tt}$ the subset of $\TET$ consisting of triples $[(\Sig,\phi,\xi)]$ where
$\phi$ is a homeomorphism which maps $\Tt$ onto an admissible triangulation of $\Sig$. Proposition \ref{ETprB1} implies that the family
$\{\U_\Tt: \; \Tt \in \mathcal{TR}(S_g,\Ah)\}$ covers the space $\TET$. We will define coordinate charts on $\U_\Tt$ for each $\Tt$ in
$\mathcal{TR}(S_g,\Ah)$.\\

\subsection{Definition of local charts}\label{ETlocch}

Given an equivalence class of triangulations $\Tt$ in $\mathcal{TR}(S_g,\Ah)$,  let $[(\Sig,\phi,\xi)]$ be a point in $\U_\Tt$. By definition,
we can assume that $\T=\phi(\Tt)$ is an admissible triangulation of $\Sig$. By Proposition \ref{ETprB2}, we know that $\T$ is
unique.\\

\noindent Slitting open  $\Sig$ along the erasing forest $\hat{A}=\phi(\Ah)$, we obtain a translation surface with piece-wise geodesic boundary
which will be denoted by $\hat{\Sig}$. The triangulation $\T$ of $\Sig$ induces a geodesic triangulation $\hat{\T}$ of $\hat{\Sig}$. Let $N_1$
be the number of edges of $\hat{\T}$, and $N_2$ be the number of triangles of $\hat{\T}$. By computing the Euler characteristic of $\Sig$, we
see that:

$$N_1=3(2g+m-2)+4(n-m) \text{ and } N_2=2(2g+m-2)+2(n-m).$$

\noindent We construct a map from $\U_\Tt$ to $\C^{N_1}$ as follows: choose an orientation for every edge of $\hat{\T}$, for each triangle
$\Delta$ in $\hat{\T}$, there exists an isometric embedding of this triangle into $\R^2$ such that the vector field $\xi$ is mapped to the
constant vertical vector field $(0,1)$ on the image of $\Delta$. By this embedding, each oriented side of the triangle $\Delta$ is mapped into a
vector in $\R^2\simeq\C$. As a consequence, we can associate to every oriented edge $e$ of $\hat{\T}$ a complex number $z(e)$. Note that, even
though each edge $e$ in the interior of $\Sig$ belongs to two distinct triangles, the complex number $z(e)$ is well defined since the vector
field
$\xi$ is parallel and normalized. The procedure above defines a map from $\U_\Tt$ into $\C^{N_1}$, we denote this map by $\Psi_\Tt$.\\

\noindent First, we have the following important observations:

\begin{lemma}\label{ETlmB2}

\begin{itemize}
\item[i)] Let $e_i, e_j, e_k$ be three edges of $\hat{\T}$ which bound a triangle. Then we have

\begin{equation}\label{TriaEq}
\pm z(e_i)\pm z(e_j)\pm z(e_k)=0,
\end{equation}

where the  signs are determined by the orientation of $e_i, e_j$ and $e_k$.\\

\item[ii)] If $e_1,\dots,e_k$ are the $k$ edges of $\hat{\T}$ which bound an open disk in $\hat{\Sig}$, then we have

\begin{equation}\label{DiskEq}
\pm z(e_1)\pm\dots\pm z(e_k)=0,
\end{equation}

where, again, the signs are determined by the orientations of the edges.\\

\end{itemize}

\end{lemma}

\dem Assertion $i)$ is straightforward. Assertion $ii)$ follows from $i)$. Namely, let $\mathrm{D}$ denote the disk bounded by $e_1,\dots,e_k$.
The disk $\mathrm{D}$ is divided into triangles by the triangulation $\hat{\T}$. By $i)$, three sides of a triangle verify (\ref{TriaEq}). Note
that every edge of $\hat{\T}$ inside $\mathrm{D}$ belongs to two distinct triangles. If, for each triangle, we choose the orientation of its
boundary coherently with the orientation of the surface, and write the corresponding equation  according to this orientation, then, by taking
the sum over all the triangles inside $\mathrm{D}$, we get (\ref{DiskEq}).\carre

\begin{lemma}\label{ETlmB3}
Let $(e,\bar{e})$ be a pair of edges in the boundary of $\hat{\Sig}$ which corresponds to an edge of a tree $A_j$ in $\hat{A}$. Suppose that $e$
and $\bar{e}$ are oriented following an orientation of $\hat{\Sig}$, then we have

\begin{equation}\label{BdrEq}
z(\bar{e})=-e^{\imath\theta}z(e)
\end{equation}

\noindent where the number $\theta$ is determined up to sign by the angles $\underline{\alpha}$, and the tree $\A_j$.
\end{lemma}

\dem Let $c$ be a path in $\hat{\Sig}$ joining the midpoint of $e$ to the midpoint of $\bar{e}$ such that
$\inter(c)\cap\partial\hat{\Sig}=\vide$. Let $\tilde{e}$ denote the edge of $A_j$ which corresponds to the pair $(e,\bar{e})$, and $q$ denote
the mid-point of $\tilde{e}$. Then $c$ corresponds to a closed curve $\gamma$ in $\Sig$ which intersects $\hat{A}$ only at $q$ transversely.
Observe that (\ref{BdrEq}) is verified when $\theta$ is the rotation angle of the holonomy of $\gamma$. Note that this angle is independent of
the choice the base-point of $\gamma$.\\

\noindent To simplify the notations, we denote by $\orth(\gamma)$ the linear part of the holonomy of $\gamma$, which is a rotation. We need to
show that the angle $\theta$ is determined up to sign by the tree $A_j$ and the angles $\alpha_1,\dots,\alpha_n$.\\

\noindent Since $A_j$ is a tree, $A_j\setm \inter(\tilde{e})$ has two connected components, let $A'_j$ denote one them, which is a sub-tree of
$A_j$. Choose an orientation of $\gamma$, and let $q_1$ and $q_2$ be two points in $\gamma$ close to $q$ so that $q$ is between $q_1$ and $q_2$.
Let $s_1$, $s_2$ denote the two sub-arcs of $\gamma$ with endpoints $q_1,q_2$, where $s_2$ is the sub-arc containing $q$. We can then find a
simple arc $s_3$ with endpoints $q_1,q_2$ such that $s_2\cup s_3$ is the boundary of a disk $\mathrm{D}$ which contains the tree $A'_j$, and no
other singular points of $\Sig$ except the vertices of $A'_j$ are contained in $\mathrm{D}$. Note that we must have

$$s_3\cap \hat{A}=\vide.$$

\noindent As a consequence, we see that $\gamma$ is (freely) homotopic to the curve $\gamma_1\cdot\gamma_2$, where $\gamma_1=s_1\cup s_3$, and
$\gamma_2=s_2\cup s_3$. Hence, the rotation angle of $\orth(\gamma)$ is equal to the sum of the rotation angles of $\orth(\gamma_1)$, and
$\orth(\gamma_2)$. By definition of erasing forest, we know that $\orth(\gamma_1)=\Id$, meanwhile, the rotation angle of $\orth(\gamma_2)$ is
equal to the sum of all the angles at the singular points contained in the disk $\mathrm{D}$ modulo $2\pi$. Therefore, we have

$$\theta=\sum_{\phi(p_i)\in A'_j}\alpha_i \mod 2\pi.$$

\noindent Observe that, the angle $\theta$ is only determined up to sign since it depends on the orientation of $\gamma$, and on the choice of
$A'_j$.\carre

Let $\STt$ denote the linear equation system consisting of $N_2$ equations of type (\ref{TriaEq}) corresponding to the triangles of $\Tt$, and
$n-m$ equations of type (\ref{BdrEq}) corresponding to the edges of $\hat{A}$. From what we have seen, the vector $\Psi_\Tt([(\Sig,\phi,\xi)])$
is a solution of the system $\STt$. Let $\VTt$ denote the subspace of $\C^{N_1}$ consisting of solutions of the system $\STt$.  We have

\begin{lemma}\label{ETlmB4}

$\Psi_\Tt(\U_\Tt)$ is an open subset of $\VTt$.

\end{lemma}

\dem The fact that $\Psi_\Tt(\U_\Tt)$ is contained in $\VTt$ is a direct consequence of Lemma \ref{ETlmB2}, and Lemma \ref{ETlmB3}.\\

\noindent Now, let $Z$ be the image of $[(\Sig,\phi,\xi)]$ by $\Psi_\Tt$, and let $Z'=(z'_1,\dots,z'_{N_1})$ be a vector in a neighborhood of
$Z$ in $\VTt$. Using the triangulation $\T$ of $\Sig$, we construct a flat surface from $Z'$ as follows:

\begin{itemize}
\item[.] Construct an Euclidean triangle from $z'_i, z'_j, z'_k$ if $z'_i,z'_j,z'_k$ verify an equation of type (\ref{TriaEq}).

\item[.] Identify two sides of two distinct triangles if they correspond to the same complex number $z'_i$.

\item[.] Identify the edges corresponding to $z_i$ and $z_j$ if $z_i$ and $z_j$ satisfy an equation of type (\ref{BdrEq}).
\end{itemize}

\noindent Clearly by this construction we obtain a translation surface $\Sig'$ homeomorphic to $\Sig$. The surface $\Sig'$ has $n$ singular
points of cone angles $\alpha_1,\dots,\alpha_n$ in the interior, there exists an erasing forest $\hat{A}'$ on $\Sig'$. Moreover, we also get a
triangulation $\T'$ of $\Sig'$ by geodesic segments. Each triangle in $\T'$ corresponds to a triangle in $\R^2$ specified by three complex
numbers which are coordinates of $Z'$, hence  we get a normalized parallel vector field $\xi'$ on $\Sig'\setm \hat{A}'$ which is defined by the
constant vertical vector field $(0,1)$ on the Euclidean plan $\R^2$.\\

\noindent Define an orientation preserving homeomorphism $\DS{f :\Sig \ra \Sig'}$ as follows: $f$ maps each edge of $\T$ onto the corresponding
edge of $\T'$ (\ie the edge of $\T$ that corresponds to the same coordinate), and the restriction $f$ on each triangle of $\T$ is a linear
transformation of $\R^2$. Let $\phi'$ denote the map $\DS{f\circ\phi}$, it follows that the triple $(\Sig',\phi',\xi')$ represents a point of
$\U_\Tt$ close to $[(\Sig,\phi,\xi)]$. By construction, it is clear that $Z'=\Psi_\Tt([(\Sig',\phi',\xi')])$. Hence, we deduce that
$\Psi_\Tt(\U_\Tt)$ is an open set of $\VTt$. \carre

\subsection{ Injectivity of $\Psi_\Tt$}

\begin{lemma}\label{ETlmB5}
The map $\Psi_\Tt$ is injective.\\
\end{lemma}

\dem Let $[(\Sig_1,\phi_1,\xi_1)]$ and $[(\Sig_2,\phi_2,\xi_2)]$ be two points in $\U_\Tt$ such that

$$\Psi_\Tt([(\Sig_1,\phi_1,\xi_1)])=\Psi_\Tt([(\Sig_2,\phi_2,\xi_2)]).$$

\noindent By definition, we can assume that $\T_1=\phi_1(\Tt)$ and $\T_2=\phi_2(\Tt)$ are admissible triangulations of $(\Sig_1,\phi_1)$ and
$(\Sig_2,\phi_2)$ respectively.\\

\noindent  Now, the hypothesis $\Psi_\Tt([(\Sig_1,\phi_1,\xi_1)])=\Psi_\Tt([(\Sig_2,\phi_2,\xi_2)])$ implies that there exists an isometry $h:
\Sig_1\ra \Sig_2$, which maps each triangle of $\T_1$ onto a triangle of $\T_2$, and also $\xi_1$ to $\xi_2$. It follows that the homeomorphism
$\phi_2^{-1}\circ h \circ \phi_1 : S_g \ra S_g$ fixes all the points in $\V$, and preserves each triangles of $\Tt$. We deduce that the map
$\phi_2^{-1}\circ h \circ \phi_1$ is isotopic to the identity of $S_g$ by an isotopy fixing all the points in $\V$. Therefore, by definition, we
have $[(\Sig_1,\phi_1,\xi_1)]=[(\Sig_2,\phi_2,\xi_2)]$.\carre

\subsection{ Computation of dimensions}

\begin{lemma}\label{ETlmB6}
$\dim_\C \VTt=\left\{%
\begin{array}{ll}
    2g+n-1, & \hbox{ if $ \alpha_i\in 2\pi\N, \; \forall i=1,\dots,n $;} \\
    2g+n-2  & \hbox{ otherwise.} \\
\end{array}%
\right.$
\end{lemma}

\dem Recall that $\VTt$ is the subspace of $\C^{N_1}$ consisting of solutions of the system $\STt$. Since the system $\STt$ contains $N_2+n-m$
equations, we have

\begin{equation}\label{ETineq1}
\dim \VTt \geq N_1-(N_2+n-m)=2g+n-2.
\end{equation}

\noindent Let $a_1,\bar{a}_1,\dots,a_{n-m}, \bar{a}_{n-m}$ denote the edges of $\hat{\T}$ which are contained in the boundary of $\hat{\Sig}$ so
that the pair $(a_i,\bar{a}_i)$ corresponds to an edge of $\hat{A}$. Choose a family of edges $\{b_1,\dots,b_{2g+m-1}\}$ in $\hat{\T}$ such that
$\DS{\inter(\hat{\Sig})\setm \cup_{j=1}^{2g+m-1}b_j}$ is an open disk.\\

\noindent  From Lemma \ref{ETlmB2} $ii)$, we deduce that if $e$ is any edge of $\hat{\T}$ which does not belong to the set

$$\{a_1,\bar{a}_1,\dots,a_{n-m},\bar{a}_{n-m},b_1,\dots,b_{2g+m-1}\},$$

\noindent then $z(e)$ can be written as a linear combination of

$$z(a_1),z(\bar{a}_1),\dots,z(a_{s_1+\dots+s_m}),z(\bar{a}_{n-m}),z(b_1),\dots,z(b_{2g+m+n-1})$$

\noindent with coefficients in $\{-1,0,1\}$. From Lemma \ref{ETlmB3}, we know that $z(\bar{a}_i)=-e^{\imath\theta_i}z(a_i)$, where $\theta_i$ is
determined by $\underline{\alpha}$ and $\Ah$. Thus, the number $z(e)$ is a linear function of

$$(z(a_1),\dots,z(a_{n-m}),z(b_1),\dots,z(b_{2g+m-1}))$$

\noindent with coefficients determined by $\underline{\alpha}$ and $\Ah$. We deduce that

\begin{equation}\label{ETineq2}
\dim \VTt \leq 2g+n-1.
\end{equation}

\noindent Suppose that the edges $a_1,\bar{a}_1,\dots,a_{n-m},\bar{a}_{n-m}$ are oriented coherently with the orientation of
$\partial\hat{\Sig}$, apply (\ref{DiskEq}) to the disk $\DS{\mathbf{D}=\inter(\hat{\Sig})\setm\cup_{j=1}^{2g+m-1}b_j}$, we get

\begin{equation}\label{ETeq1}
\sum_{i=1}^{n-m} (z(a_i)+z(\bar{a}_i))=\sum_{i=1}^{n-m}(1-e^{\imath\theta_i})z(a_i)=0
\end{equation}

\noindent The numbers $z(b_j), \; j=1,\dots,2g+m-1,$ do not appear in the equation (\ref{ETeq1}) because each of the edges $b_j$ belongs to two
distinct triangles. Here, we have two issues:

\begin{itemize}
\item[$\bullet$] Case 1: there exists $i\in\{1,\dots,n\}$ such that $\alpha_i\notin 2\pi\N$. The equation (\ref{ETeq1}) is then non-trivial,
which means that the vector $(z(a_1),\dots,z(a_{n-m}),z(b_1),\dots,z(b_{2g+m-1}))$ belongs to a hyperplane of $\C^{2g+n-1}$. Therefore we have

\begin{equation}\label{ETineq3}
\dim \VTt\leq 2g+n-2.
\end{equation}

From (\ref{ETineq1}) and (\ref{ETineq3}), we conclude that $\dim_\C \VTt=2g+n-2$.\\

\item[$\bullet$] Case 2: $\alpha_i \in 2\pi\N$ for every $i$ in $\{1,\dots,n\}$. In this case, the equation (\ref{ETeq1}) is trivial. However,
this also means that the sum of all equations in the system $\STt$, with appropriate choices of signs, is the trivial equation $0=0$. This
implies $\mathrm{rk}(\STt)\leq N_2+(n-m)-1$. Hence

\begin{equation}\label{ETineq4}
    \dim \VTt\geq N_1-(N_2+n-m-1)=2g+n-1.
\end{equation}

From (\ref{ETineq2}) and (\ref{ETineq4}), we conclude that $\dim \VTt=2g+n-1$.\\
\end{itemize}

The lemma is then proved.\carre

\subsection{ Coordinate change}

Let $\Tt_1, \Tt_2$ be two equivalence classes of triangulations in $\mathcal{TR}(S_g,\Ah)$. Suppose that $\U_{\Tt_1}\cap \U_{\Tt_2}\neq \vide$,
and let $[(\Sig,\phi,\xi)]$ be a point in $\U_{\Tt_1}\cap \U_{\Tt_2}\neq \vide$. Let $\T_1,\T_2$ be the admissible triangulations of $\Sig$
corresponding to $\Tt_1$ and $\Tt_2$ respectively. As usual, we denote by $\Psi_{\Tt_1},\Psi_{\Tt_2}$ the local charts on $\U_{\Tt_1}$ and
$\U_{\Tt_2}$ respectively. We have:

\begin{lemma}\label{lmA10}
There exists an invertible complex linear map

$$\mathbf{L}: \C^{N_1}\lra \C^{N_1}$$

\noindent such that $\Psi_{\Tt_2}([(\Sig',\phi',\xi')])=\mathbf{L}\circ\Psi_{\Tt_1}([(\Sig',\phi',\xi')])$, for every $[(\Sig',\phi',\xi')]$
in a neighborhood of $[(\Sig,\phi,\xi)]$.\\
\end{lemma}

\dem Let $e$ be an edge of $\T_2$. Let $\{\Delta_i, \; i\in I\}$ be the set of triangles in $\T_1$ which are crossed by $e$, that is
$\Delta_i\cap\inter(e)\neq\vide, \; \forall i\in I$. Using the developing map, we can construct a  polygon $\P_e$ in $\R^2$ by gluing
successively isometric copies of $\Delta_i$'s ($i\in I$) so that $e$ corresponds to a diagonal $\tilde{e}$ inside $\P_e$. Note that a triangle
$\Delta_i$ may have several copies in the polygon $\P_e$. We will call $\P_e$ the {\em developing polygon of $e$ with respect to $\T_1$}.\\

\noindent  By this construction, we also get a map $\DS{\varphi: \P_e\lra \Sig}$ which is locally isometric, such that $\varphi(\tilde{e})=e$.
The inverse images of the edges of $\T_1$ under $\varphi$ give rise to a triangulation of $\P_e$ by diagonals. Note that by construction, all of
the diagonals of this triangulation intersect the diagonal $\tilde{e}$.\\

\noindent Since the map $\varphi$ sends segments in the boundary of $\P_e$ onto edges of $\T_1$, it follows that the complex numbers associated
to the edge $e$ can be written as linear function of the complex numbers associated to some edges of $\T_1$. The coefficients of these linear
functions are unchanged if we replace $[(\Sig,\phi,\xi)]$ by another point $[(\Sig',\phi',\xi')]$ nearby in $\U_{\Tt_1}\cap \U_{\Tt_2}$, and
this argument is reciprocal between $\T_1$ and $\T_2$. We deduce that the coordinate change between $\Psi_{\Tt_1}$ and $\Psi_{\Tt_2}$, in a
neighborhood of $[(\Sig,\phi,\xi)]$, is a complex linear transformation of $\C^{N_1}$ which sends $\mathrm{V}_{\Tt_1}$ onto
$\mathrm{V}_{\Tt_2}$. The lemma is then proved. \carre


\section{Action of Mapping Class Group}\label{ETMCGASec}

In this section, we will prove that the action of $\GgA$ on $\TET$ is properly discontinuous. For this purpose, we first define a map $\Xi$ from
$\TET$ into $\Tgn$ as follows: the image of a point $[(\Sig,\phi,\xi)]$ in $\TET$ under $\Xi$ is the point in $\Tgn$ represented by the pair
$(\Sig,\phi)$, where $\Sig$ is now considered as a Riemann surface with $n$ marked points $\{\phi(p_1),\dots,\phi(p_n)\}$. \\

\noindent Observe that the group $\GgA$ is naturally embedded into the mapping class group $\Ggn$ of the punctured surface
$S_g\setm\{p_1,\dots,p_n\}$, since by definition, we have

$$\HomgAo=\HomgA \cap \Homgno,$$

\noindent where $\Homgno$ is the set of orientation preserving homeomorphisms of $S_g$ which are isotopic to $\Id_{S_g}$ relative to
$\{p_1,\dots,p_n\}$. Clearly, the map $\Xi$ is equivariant with respect to the actions of $\GgA$ on $\TET$, and on $\Tgn$. It is well known that
the group $\Ggn$ acts properly discontinuously on $\Tgn$, therefore, it suffices to prove

\begin{proposition}\label{ETprB3}
The map $\DS{\Xi: \TET \lra \Tgn}$ is continuous.
\end{proposition}

\dem Let $[(\Sig,\phi,\xi)]$ be a point in $\TET$, and $\{[(\Sig_k,\phi_k,\xi_k)], \; k\in \N\}$ be a sequence in $\TET$ converging to
$[(\Sig,\phi,\xi)]$.  The point $[(\Sig,\phi,\xi)]$ belongs to  a set $\U_\Tt$, where $\Tt$ is an equivalence class in $\mathcal{TR}(S_g,\Ah)$.
Without loss of generality, we can assume that the sequence $\{[(\Sig_k,\phi_k,\xi_k)], \; k\in \N \}$ is also contained in $\U_\Tt$. Let
$\Psi_\Tt$ be the local chart of $\TET$ which is defined on $\U_\Tt$. Put $\DS{Z =\Psi_\Tt([(\Sig,\phi,\xi)])}$, and
$\DS{Z_{k}=\Psi_\Tt([(\Sig_k,\phi_k,\xi_k)])}$, by assumption we have $\DS{Z_k\overset{k\rightarrow\infty}{\lra} Z}$ in $\C^{N_1}$.\\

\noindent Let $\T$ be the admissible triangulation of $\Sig$ corresponding to $\Tt$. Recall that, by the definition of $\Psi_\Tt$, for every
point $[(\Sig',\phi',\xi')]$ in $\U_\Tt$, we can write $\phi'= f\circ\phi$, where $f: \Sig\lra \Sig'$ is a homeomorphism verifying

\begin{itemize}
\item[$\bullet$] $f$ sends each edge of $\T$ onto an edge of an admissible triangulation $\T'$ of $\Sig'$,

\item[$\bullet$] the restriction of $f'$ into the a triangle of $\T$ is a linear transformation of $\R^2$.

\end{itemize}

\noindent Therefore, for every $k\in \N$, we can assume that $\phi_k=f_k\circ\phi$, where $f_k: \Sig \lra \Sig_k$ is a homeomorphism with the
properties above. It is clear that, as $Z_k$ tends to $Z$, the restriction of $f_k$ on each triangle of $\T$ tends to identity, which implies
$\DS{\lim_{k\rightarrow\infty} K(f_k)=1}$, where $K(f_k)$ is the dilatation of $f_k$. By the definition of Teichmüller distance
$\di_{\mathrm{Teich}}$, it follows that

$$\lim_{k\rightarrow\infty}\di_{\mathrm{Teich}}([(\Sig,\phi)],[(\Sig_k,\phi_k)])=0,$$

\noindent and the proposition follows.\carre

Thank to a result of Troyanov \cite{Troy}, more can be said about the map $\Xi$. Observe that the map $\Xi$ can be defined on $\TETv$, and we
have

\begin{proposition}\label{ETprB4}
The restriction of $\Xi$ to $\TETvi$ is injective.
\end{proposition}

\dem Let $[(\Sig_1,\phi_1)]$ and $[(\Sig_2,\phi_2)]$ be two points in $\TETvi$ such that $(\Sig_1\setm\{\phi_1(p_1),\dots,\phi_1(p_n)\},\phi_1)$
and $(\Sig_2\setm\{\phi_2(p_1),\dots,\phi_2(p_n)\},\phi_2)$ belong to the same equivalence class in $\Tgn$, we have to prove that
$\DS{[(\Sig_1,\phi_1)]=[(\Sig_2,\phi_2)]}$ in $\TETv$.\\

\noindent By assumption, there exists a conformal homeomorphism $h:\Sig_1\longrightarrow \Sig_2$ such that $\DS{\phi_2^{-1}\circ h \circ
\phi_1}$ is isotopic to $\Id_{S_g}$ by an isotopy which is identity on the set $\{p_1,\dots,p_n\}$. First, we prove that $h$ is also an isometry
between the two flat surfaces $\Sig_1$ and $\Sig_2$.\\

\noindent Consider $\Sig_1$ and $\Sig_2$ as Riemann surfaces equipped with conformal flat metrics $f_1$ and $f_2$ respectively. Let
$\mathbf{div}_i$ denote the formal sum $\DS{\sum_{k=1}^n s_k\phi_i(p_k)}$, where $\DS{s_k=\frac{\alpha_k}{2\pi}-1}$, we say that the metric
$f_i$ represents the {\em divisor} $\mathbf{div}_i, \; i=1,2$.\\

\noindent Since $h$ is a conformal homeomorphism, it follows that $h^*f_2$ is also a conformal flat metric on $\Sig_1$. By assumption, we have
$h(\mathbf{div}_1)=\mathbf{div}_2$, and we deduce that $h^*f_2$ represents $\mathbf{div}_1$ too.  Now, from Proposition 2 \cite{Troy}, there
exists $\lambda>0$ such that $f_1=\lambda\ h^*f_2$. Since we have $\mathrm{Area}_{f_1}(\Sig_1)=\mathrm{Area}_{f_2}(\Sig_2)=1$, it follows that
$\lambda=1$. Therefore we have $f_1=h^*f_2$, in other words, $h$ is an isometry from the flat surface $\Sig_1$ onto the flat surface $\Sig_2$.\\

All we need  to prove now is that $\phi_2^{-1}\circ h \circ \phi_1$ preserves the forest $\Ah$. Since $h$ is an isometry of flat surfaces,
$h(\phi_1(\hat{\A}))$ is a union of geodesic trees whose vertices are singular points of $\Sig_2$. Let $a$ be an edge of a tree in $\Ah$, then
$\phi_1(a)$ is a geodesic segment on $\Sig_1$, hence $h(\phi_1(a))$ is a geodesic segment of $\Sig_2$. By definition, $\phi_2(a)$ is also a
geodesic segment of $\Sig_2$ which has the same endpoints as $h(\phi_1(a))$.\\

\noindent  By assumption, there exists an isotopy relative to $\{p_1,\dots,p_n\}$ from $h\circ \phi_1$ to $\phi_2$, it follows that $\phi_2(a)$
and $h(\phi_1(a))$ are homotopic with fixed endpoints in $\Sig_2$ relative to $\{\phi_2(p_1),\dots,\phi_2(p_n)\}$. Now, from Lemma \ref{ETlmA},
we have $h(\phi_1(a))=\phi_2(a)$. Since this is true for every edges in $\hat{\A}$, we conclude that $h\circ\phi_1(\hat{\A})=\phi_2(\hat{\A})$,
or equivalently, $\phi_2^{-1}\circ h \circ \phi_1 (\hat{\A})=\hat{\A}$. It follows immediately that $\phi_2^{-1}\circ h \circ \phi_1 \in
\HomgAo$, in other words, $(\Sig_1,\phi_1)$ and $(\Sig_2,\phi_2)$ are equivalent in $\TETvi$.\carre

Since $\TET$ is a $\C^*$-bundle over $\TETvi$, we can then consider $\TET$ as a $\C^*$-bundle over the subset $\Xi(\TETvi)$ of $\Tgn$.\\


\section{Changes of triangulations}\label{TriaCh}

Let $[(\Sig,\phi)]$ be an element of the space $\TETv$, we have seen that an admissible triangulation of $(\Sig,\phi)$ (cf. Definition
\ref{defAT}) allows us to construct a local chart for $\TET$. In this section, we are interested in relations between geodesic triangulations of
$\Sig$. More precisely, we want to answer the question: How to go from an admissible triangulation to another one. This will play a crucial role
in our construction of the volume form on $\TET$.\\

Let us start with a simple example: let $ABCD$ be a convex quadrilateral in $\R^2$, there are only two ways to triangulate $ABCD$: one by
adding the diagonal $AC$, and the other by adding the diagonal $BD$.\\

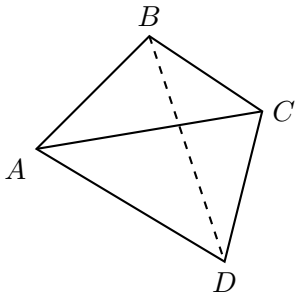
\begin{center}
\begin{tikzpicture}[scale=0.5]
\draw[black, thick] (0,0) -- (3,3) -- (6,1) -- (5,-3) -- cycle; \draw[black,thick] (0,0) -- (6,1); \draw (0,0) node[below left] {$A$}; \draw
(3,3) node[above] {$B$}; \draw (6,1) node[right] {$C$}; \draw (5,-3) node[below] {$D$}; \draw[thick, dashed] (3,3) -- (5,-3);

\end{tikzpicture}
\end{center}

\noindent This example suggests

\begin{definition}[Elementary Move and Related Triangulations]
\label{defEM} Let $\Sig$ be a flat surface with geodesic boundary. Let $\T$ be a triangulation of $\Sig$ by geodesic segments whose set of
vertices contains the set of singularities of $\Sig$. An {\em elementary move} of $\T$ is a transformation as follows: take two adjacent
triangles of $\T$ which form a convex quadrilateral, replace the common side of the two triangles by the other diagonal of the quadrilateral (if
these two triangles have more than one common side, just take one of them). After such a move, we obtain evidently a another geodesic
triangulation of $\Sig$ with the same set of vertices as $\T$.\\

\noindent Let $\T_1,\T_2$ be two geodesic triangulations of $\Sig$ whose sets of vertices coincide. We say that $\T_1$ and $\T_2$ are {\em
related} if there exists a sequence of elementary moves which transform $\T_1$ into $\T_2$.\\

\end{definition}

The main result of this section is the following theorem

\begin{theorem}\label{ETthF}
Let $\Sig$ be a flat surface with piece-wise geodesic boundary. Let $V$ be a finite subset of $\Sig$ which contains all the singularities.
Suppose that $\Sig$ satisfies the following condition

$$\pQ \hspace{1cm} \text{ for every closed curve } c \subset \inter(\Sig)\setm V,\text{ we have } \orth(c)\in\{\pm\Id\},$$

\noindent where $\orth(c)$  is the orthogonal part of the holonomy of $c$. Let $\T_1, \T_2$ be two geodesic triangulations of $\Sig$ such that
the set of vertices of $\T_i$ is $V$, $i=1,2$, then $\T_1$ and $\T_2$ are related.\\
\end{theorem}

\rem The changes of triangulations by elementary moves, which are also called {\em flips}, are already studied in the context of general flat
surfaces (not necessarily translation surfaces). In the general situation, Theorem \ref{ETthF} is already known, it results from the fact that
any geodesic triangulation whose vertex set contains all the singularities can be transformed by flips into a special one, called Delaunay
triangulation, which is unique up to some flips (see \cite{BobSpr} for further detail).\\

\noindent For the case at hand, we have an elementary proof of this fact which is based on an observation on polygons, and uses some basic
properties of translation and half-translation surfaces, this proof is given in Appendices, Section \ref{prfFSect}.\\


\section{Volume form}\label{ETVolSec}

In this section, we define a volume form on the manifold $\TET$. Recall that, if $\mathrm{L}: \mathrm{E}\ra \mathrm{F}$ is a linear map between
(real) vector spaces which is surjective, then given a volume form $\mu_\mathrm{E}$ on $\mathrm{E}$, and a volume form $\mu_\mathrm{F}$ on
$\mathrm{F}$, one can define a volume form $\mu$ on $\ker(\mathrm{L})$ as follows: let $\iota$ be the embedding of $\ker(\mathrm{L})$ into
$\mathrm{E}$, the volume form $\mu$ on $\ker(\mathrm{L})$ is defined to be $\iota^*\tilde{\mu}$, where $\tilde{\mu}$ is any element of
$\Lambda^{(\dim\mathrm{E}-\dim\mathrm{F})}\mathrm{E}^*$ such that:

$$ \mu_\mathrm{E}=\tilde{\mu}\wedge\mathrm{L}^*\mu_\mathrm{F}.$$


\subsection{Definitions}

Let $\Tt$ be a triangulation of $S_g$, which represents an equivalence class in $\mathcal{TR}(S_g,\Ah)$. Cut open $S_g$ along the edges of the
forest $\Ah$, and denote the new surface by $\hat{S}_g$. Let $\hat{\Tt}$ denote the triangulation of $\hat{S}_g$ which is induced by $\Tt$. As
usual, let $N_1,N_2$ denote the number of edges, and the number of triangles in $\hat{\Tt}$ respectively. Let $\Psi_\Tt : \U_\Tt\lra \C^{N_1}$
be the local chart associated to $\Tt$. Recall that $\Psi_\Tt(\U_\Tt)$ is an open subset of the solution space $\VTt$ of a system $\STt$, which
consists of $N_2$ equations of type (\ref{TriaEq}) and $(n-m)$ equations of type (\ref{BdrEq}).\\

\noindent Let $a_1,\dots,a_{N_2+(n-m)}$ denote the vectors of $(\C^{N_1})^*$ which correspond to the equations of the system $\STt$. A vector
$a_i$ is said to be {\em normalized} if each of its coordinates is either $0$, or a complex number of module $1$. Consider the complex linear
map $\ATt: \C^{N_1}\lra \C^{N_2+(n-m)}$, which is defined in the canonical basis of $\C^{N_1}$ and $\C^{N_2+(n-m)}$ by the matrix

$$\ATt=\left(%
\begin{array}{c}
  a_1 \\
  \vdots\\
  a_{N_2+(n-m)} \\
\end{array}%
\right).$$

\noindent We have $\VTt=\ker\ATt$. The map $\ATt$ is said to be {\em normalized} if each row of its matrix in the canonical basis is normalized.
We have two cases:

\begin{itemize}
\item[$\bullet$]\underline{Case 1:} there exist $i\in\{1,\dots,n\}$ such that $\alpha_i\notin 2\pi\N$. In this case, we have seen that (cf.
Lemma \ref{ETlmB6}) $\dim \VTt=2g+n-2$, hence, $\mathrm{rk}(\ATt)=\mathrm{rk}(\STt)=N_2+(n-m)$. Let $\lambda_{2N_1}$ et $\lambda_{2(N_2+(n-m))}$
denote the Lebesgue measures on $\C^{N_1}\simeq \R^{2N_1}$ and $\C^{N_2+(n-m)}\simeq \R^{2(N_2+(n-m))}$ respectively. Since $\ATt$ is
surjective, $\lambda_{2N_1}$ and $\lambda_{2N_2}$ induce a volume form $\nu_{\Tt}$ on $\VTt$ via the following exact sequence:

\begin{equation}\label{ExSeq1}
0\lra \VTt \hookrightarrow \C^{N_1}\overset{\ATt}{\lra}\C^{N_2+(n-m)}\lra 0
\end{equation}

\item[$\bullet$] \underline{Case 2:} for every  $i\in\{1,\dots,n\}, \; \alpha_i \in 2\pi\N$. In this case, $\dim \VTt=2g+n-1$, and
$\mathrm{rk}(\ATt)=\mathrm{rk}(\STt)=N_2+(n-m)-1$. If the vectors $a_1,\dots,a_{N_2+(n-m)}$ are normalized, and if their signs are chosen
suitably, we have $a_1+\dots+a_{N_2}=0$. Thus, without loss of generality, we can assume that $\mathrm{Im} \ATt=\mathbf{W}$, where $\mathbf{W}$
is the complex hyperplane of $\C^{N_2+(n-m)}$ defined by

$$\mathbf{W}=\{(z_1,\dots,z_{N_2+(n-m)})\in \C^{N_2+(n-m)} \; | \; z_1+\dots+z_{N_2+(n-m)}=0\}.$$

Let $\tilde{\lambda}_{2(N_2+(n-m)-1)}$ denote the volume form of $\mathbf{W}$ which is induced by the Lebesgue measure of $\C^{N_2+(n-m)}$. The
volume forms $\lambda_{2N_1}$ and $\tilde{\lambda}_{2(N_2+(n-m)-1)}$ induce a volume form $\nu_{\Tt}$ on $\VTt^*$ via the exact sequence:

\begin{equation}\label{ExSeq2}
0\lra \VTt \hookrightarrow \C^{N_1}\overset{\ATt}{\lra} \mathbf{W}\lra 0
\end{equation}

In other words, in this case, $\nu_{\Tt}$ is defined by the torsion of the following exact sequence:

\begin{equation}\label{ExSeq2b}
0\lra \VTt \hookrightarrow \C^{N_1}\overset{\ATt}{\lra} \C^{N_2+(n-m)}\overset{\mathbf{s}}{\lra}\C \lra 0
\end{equation}

\noindent where $\mathbf{s}((z_1,\dots,z_{N_2+(n-m)}))=z_1+\dots+z_{N_2+(n-m)}$, together with the Lebesgue measures on
$\C^{N_1},\C^{N_2+(n-m)},$ and $\C$.\\

\end{itemize}

In both cases, let $\mu_{\Tt}$ denote the volume form $\Psi_{\Tt}^*\nu_{\Tt}$ which is defined on $\U_\Tt$.\\

\subsection{Invariance by coordinate changes}

Let $\Tt_1$, and $\Tt_2$ be two triangulations of $S_g$ which represent two different equivalence classes in $\mathcal{TR}(S_g,\Ah)$. Assume
that $\U_{\Tt_1}\cap\U_{\Tt_2}\neq \vide$. Then we have

\begin{lemma}\label{ETlmB7}

$\mu_{\Tt_1}=\mu_{\Tt_2}$ on $\U_{\Tt_1}\cap\U_{\Tt_2}$.

\end{lemma}

\dem Let $[(\Sig,\phi,\xi)]$ be a point in $\U_{\Tt_1}\cap\U_{\Tt_2}$, and let $\T_1$, $\T_2$ be the admissible triangulations of $(\Sig,\phi)$
corresponding to $\Tt_1$ and $\Tt_2$ respectively. By Theorem \ref{ETthF}, we can assume that $\T_2$ is obtained from $\T_1$ by only one
elementary move. Let $Z_i=\Psi_{\Tt_i}([(\Sig,\phi,\xi)]),$ it is clear that the coordinates of $Z_2$ are linear functions of $Z_1$ with integer
coefficients and vice versa. We deduce that coordinates change between $\Psi_{\Tt_1}$, and $\Psi_{\Tt_2}$ is realized by (complex) linear
transformation $\mathbf{F}$ of $\C^{N_1}$, that is to say we have the following commutative  diagram:

$$\begin{array}{ccccccccc}
  0 & \lra & \mathrm{V}_{\Tt_1} & \hookrightarrow & \C^{N_1} & \overset{\mathbf{A}_{\Tt_1}}{\lra} & \mathbf{X} & \lra & 0 \\
   &  & \downarrow \mathbf{H} &  & \downarrow \mathbf{F} &  & \|\Id &  &  \\
  0 & \lra & \mathrm{V}_{\Tt_2} & \hookrightarrow & \C^{N_1} &\overset{\mathbf{A}_{\Tt_2}}{\lra} & \mathbf{X} & \lra & 0\\
\end{array}$$

\noindent where $\mathbf{X}$ is either $\C^{N_2+(n-m)}$, or $\mathbf{W}$, and $\mathbf{H}$ is the coordinates change between $\Psi_{\Tt_1}$, and
$\Psi_{\Tt_2}$. Note that the map $\mathbf{F}$ is written in the canonical basis of $\C^{N_1}$ as an integer matrix, and so is its inverse. As a
matter of fact, we can arrange so that

$$\mathbf{A}_{\Tt_1}=\left(%
\begin{array}{ccccccc}
  -1 & 1 & 1 & 0 & 0 & \cdots & 0 \\
  1 & 0 & 0 & 1 & 1 & \cdots & 0 \\
  0 & * & * & * & * & \cdots & * \\
  \cdots & \cdots & \cdots & \cdots & \cdots & \cdots & \cdots \\
  0 & * & * & * & * & \cdots & * \\
\end{array}
\right); \; %
\mathbf{A}_{\Tt_2}=\left(%
\begin{array}{ccccccc}
  -1 & 0 & 1 & 1 & 0 & \cdots & 0 \\
  1 & 1 & 0 & 0 & 1 & \cdots & 0 \\
  0 & * & * & * & * & \cdots & * \\
  \cdots & \cdots & \cdots & \cdots & \cdots & \cdots & \cdots \\
  0 & * & * & * & * & \cdots & * \\
\end{array},%
\right).$$

\noindent and

$$\mathbf{F}=\left(%
\begin{array}{ccccccc}
  1 & -1 & 0 & 1 & 0 & \cdots & 0 \\
  0 & 1 & 0 & 0 & 0 & \cdots & 0 \\
  0 & 0 & 1 & 0 & 0 & \cdots & 0 \\
  \cdots & \cdots & \cdots & \cdots & \cdots & \cdots & \cdots \\
  0 & 0 & 0 & 0 & 0 & \cdots & 1 \\
\end{array}%
\right).$$

\noindent in the canonical basis of $\C^{N_1}$ and $\C^{N_2+(n-m)}$, where for $j \geq 3$, the $j$-th rows of $\mathbf{A}_{\Tt_1}$, and
$\mathbf{A}_{\Tt_2}$ are the same.  Note that the first two rows of $\mathbf{A}_{\Tt_1}$, and $\mathbf{A}_{\Tt_2}$ correspond to two triangles
in the quadrilateral where the elementary move occurs. Since $|\det\mathbf{F}|=1$, it follows that $\DS{\nu_{\Tt_1}=\mathbf{H}^*\nu_{\Tt_2}}$,
and the lemma follows.\carre


\subsection{Invariance by action of $\GgA$}

Lemma \ref{ETlmB7} implies that the volume forms $\{\mu_\Tt : \; \Tt\in \mathcal{TR}(S_g,\Ah)\}$ give a well defined volume form, which will be
denoted by $\mu_\mathrm{Tr}$, on $\TET$. To complete the proof of Theorem \ref{ETthB}, we need the following

\begin{lemma}\label{ETlmB8}
The volume form $\mu_\mathrm{Tr}$ is invariant by the action of $\GgA$.\\
\end{lemma}

\dem This lemma follows readily from the arguments of Lemma \ref{ETlmB7}.\carre

\noindent The proof of Theorem \ref{ETthB} is now complete.\carre

\section{Proof of Proposition \ref{ETprC}}\label{prfCSec}

Let $(M,\omega)$ be a pair in $\Hg$, and let $\Sig$ denote the induced translation surface. Let $x_1,\dots,x_n$ denote the singularities of
$\Sig$ so that the cone angle at $x_i$ is $2\pi(k_i+1)$. The vertical geodesic flow determined by $\omega$ is induced by a unitary parallel
vector field $\xi$ on $\Sig\setm\{x_1,\dots,x_n\}$. The pair $(M,\omega)$ in $\Hg$ is then identified to the point $(\Sig,
\{x_1,\dots,x_n\},\xi)$ in $\MET$, where $\Ah$ is the union of $n$ points.\\

\noindent Let $\T$ be a geodesic triangulation of $\Sig$ whose vertex set is $\{x_1,\dots,x_n\}$. Note that, in this case, any geodesic
triangulation whose set of vertices coincides with the set of singularities is admissible. We  call a set of $2g+n-1$ edges of $\T$ such that
the complement of the union of those edges is a topological open disk a {\em family of primitive edges} of $\T$. Remark  that such a family
always exists since it corresponds to the complement of a maximal tree, \ie a tree which contains all the vertices, in the dual graph of $\T$.
Let $\{b_1,\dots,b_{2g+n-1}\}$ be a family of primitive edges of $\T$. Observe that $\{b_1,\dots,b_{2g+n-1}\}$ is a basis of the group
$H_1(\Sig,\{x_1,\dots,x_n\},\Z)$.\\

\noindent Let $\phi: S_g \ra \Sig$ be a homeomorphism which maps $p_i$ to $x_i, \; i=1,\dots,n$, and let $\Tt$ denote triangulation
$\phi^{-1}(\T)$ of $S_g$ whose vertex set is $\{p_1,\dots,p_n\}$. Let $\Psi_\Tt$ be the local chart associated to $\Tt$. As usual, let $\STt$
denote the system of linear equations associated to $\Tt$, $\VTt$ denote the space of solutions of $\STt$, and $\ATt$ denote the normalized
linear map associated to $\STt$. We can assume that

$$\mathrm{Im}\ATt=\mathbf{W}=\{(z_1,\dots,z_{N_2}) \in \C^{N_2} | \; z_1+\dots+z_{N_2}=0 \}.$$

\noindent Recall that in this case, $\dim_\C \VTt=2g+n-1$. Under $\Psi_\Tt$, a neighborhood of $(\Sig,\{x_1,\dots,x_n\},\xi)$ in $\MET$ is
identified to an open subset of $\VTt$. There exists a neighborhood $\U$ of $(\Sig,\{x_1,\dots,x_n\},\xi)$ such that, for any point
$(\Sig',\{x'_1,\dots,x'_n\},\xi')$ in $\U$, there exists a homeomorphism $f_{\Sig'}: \Sig\lra \Sig'$ such that $f_{\Sig'}(\T)=\T'$ is an
admissible triangulation of $\Sig'$. Let $b'_i$ denote $f_{\Sig'}(b_i) , \; i=1,\dots,2g+n-1$, then the segments $\{b'_1,\dots,b'_{2g+n-1}\}$
form a basis of the group $H_1(\Sig',\{x'_1,\dots,x'_n\},\Z)$. Hence, we can define a local chart of $\Hg$ by the following period mapping

$$\begin{array}{cccc}
\Phi : & \U & \lra &\C^{2g+n-1}\\
   & (\Sigma',\{x'_1,\dots,x'_n\},\xi')\simeq (M',\omega') & \longmapsto & (\int_{b'_1}\omega',\dots,\int_{b'_{2g+n-1}}\omega')\\

\end{array}$$

\noindent By the construction of $\Psi_\Tt$, we can assume that, if $\Psi_\Tt(\Sig',\{x'_1,\dots,x'_n\},\xi')=(z_1,\dots,z_{N_1})$, then the
complex numbers $z_1,\dots,z_{2g+n-1}$ are associated to the edges $b'_1,\dots,b'_{2g+n-1}$. It follows that the map $\DS{\Psi_\Tt\circ
\Phi^{-1}: \Phi(\U)\subset \C^{2g+n-1}\lra\C^{N_1}}$ is a restriction of a linear map from $\C^{2g+n-1}$ into $\C^{N_1}$ which is injective,
hence, $\Psi_\Tt\circ \Phi^{-1}$ is a restriction to $\Phi(\U)$ of an isomorphism from $\C^{2g+n-1}$ onto $\VTt$.\\

\noindent By definition, $\mu_0=\Phi^*\lambda_{2(2g+n-1)}$, where $\lambda_{2(2g+n-1)}$ is the Lebesgue measure of $\C^{2g+n-1}$, and
$\mu_{\mathrm{Tr}}=\Psi_\Tt^*\nu_\Tt$, where $\nu_\Tt$ is  the volume form on $\VTt$ which is defined by the exact sequence (\ref{ExSeq2}).
Clearly, on $\C^{2g+n-1}$ we have

$$(\Psi_\Tt\circ\Phi^{-1})^*\nu_\Tt=\lambda\lambda_{2g+n-1},$$

\noindent where $\lambda$ is a non-zero constant. This implies $\mu_{\mathrm{Tr}}=\lambda\mu_0$ on a neighborhood of
$(\Sig,\{x_1,\dots,x_n\},\xi)$. We deduce that $\displaystyle{ \mu_{\mathrm{Tr}}\over \mu_0}$ is locally constant, consequently,
$\displaystyle{\mu_{\mathrm{Tr}}\over \mu_0}$ is constant on every connected component of $\Hg$. \carre

\section{Flat complex affine structure on moduli space of flat surfaces of genus zero}\label{SSAffStSec}

\subsection{Existence of erasing forest}
To see that spherical flat surfaces are locally a special case of flat surfaces with erasing forest, let us prove the following

\begin{proposition}
\label{ETprEx}

Let $\Sig$ be compact flat surface without boundary. Let $\{p_1,\dots,p_n\}$ denote the  singularities of $\Sig$, then there exists a geodesic
tree whose vertex set is $\{p_1,\dots,p_n\}$.

\end{proposition}

\dem Let  $C_1$ be a path of minimal length from $p_1$ to $p_2$. The path $C_1$ is a finite union of geodesic segments whose endpoints are
singular points of $\Sig$. Apart from $p_1$ and $p_2$,  $C_1$ can contain other points in $\{p_1,\dots,p_n\}$. Since $C_1$ is a path of minimal
length, it has no self intersections. By renumbering the set of singular points if necessary, we can assume that $C_1$ is a path joining $p_1$
to $p_r$, and contains points $p_2,\dots,p_{r-1}$. Note that for every point $p\in C_1$, the length of the path  from $p_1$ to $p$ along $C_1$
is the distance $\di(p_1,p)$ between them, where $\di$ is the distance on $\Sig$ which is induced by the flat metric.\\

\noindent  If $r=n$, then we have obtained a geodesic tree whose vertices are $\{p_1,\dots,p_n\}$. Assume that $r<n$, let $C_2$ be a path of
minimal length from $p_1$ to $p_{r+1}$. If $C_1 \cap C_2 =\{p_1\}$, then we get a geodesic tree whose vertex set contains at least $r+1$ points
in $\{p_1,\dots,p_n\}$. If this is not the case, let us prove that $C_2$ can not intersect $C_1$ transversely at a regular point.

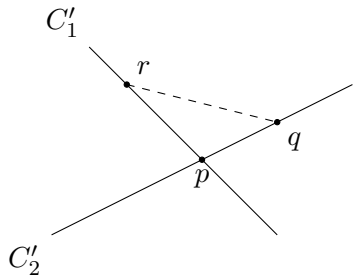
\begin{center}
\begin{tikzpicture}[scale=0.5]
\draw (-3,3) -- (2,-2); \draw (-4,-2) -- (4,2); \filldraw [black] (0,0) circle (2pt) (-2,2) circle (2pt) (2,1) circle (2pt); \draw[dashed]
(-2,2) -- (2,1);

\draw (0,0) node[anchor=north] {$p$}; \draw (-2,2) node[above right] {$r$}; \draw (2,1) node[below right] {$q$}; \draw (-3,3) node[above left]
{$C'_1$}; \draw (-4,-2) node[below left] {$C'_2$};
\end{tikzpicture}
\end{center}

\noindent Suppose that $p$ is a regular point where $C_2$ intersects $C_1$ transversely. Let $V$ be a neighborhood of $p$ such that $S_1=V\cap
C_1$ and $S_2=V\cap C_2$ are two geodesic segments, and $p$ is the unique common point of $S_1$ and $S_2$. Let $C'_1$ be the paths from $p_1$ to
$p$ along $C_1$ and $C'_2$ be the path from $p_1$ to $p$ along $C_2$, we have $\DS{\leng(C'_1)=\leng(C'_2)=\di(p_1,p)}$.\\

\noindent Let $q$ be a point in $S_2\setm C'_2$, and $r$ be a point in $S_1\cap C'_1$. Let $\overline{pq}$ denote the sub-segment of $S_2$ whose
endpoints are $p$ and $q$, and $\overline{pr}$ denote the sub-segments of $S_1$ whose endpoints are $p$ and $r$. We have
$\DS{\di(p_1,q)=\di(p_1,p)+\leng(\overline{pq})}$, and $\DS{\di(p_1,p)=\di(p_1,r)+\leng(\overline{pr})}$.\\

\noindent Since $p$ is a regular point of $\Sig$, if we choose the points $q$ and $r$ close enough to $p$, the geodesic segment $\overline{qr}$
joining $q$ and $r$ will be contained in the neighborhood $V$, and we have
$\DS{\leng(\overline{qr})<\leng(\overline{pr})+\leng(\overline{pq})}$. It follows that

$$\di(p_1,q)=\di(p_1,r)+\leng(\overline{pr})+\leng(\overline{pq})> \di(p_1,r)+\leng(\overline{qr}).$$

\noindent The above inequality is in contradiction with the definition of the distance $\di$. Thus,  we conclude that $C_2$ cannot intersect
$C_1$ transversely  at a regular point. This implies that the last intersection point of $C_1$ and $C_2$, that is the intersection point of
furthest distance from $p_1$, must be a singular point $p_k$ of $\Sig$. Omit the part of $C_2$ from $p_1$ to $p_k$ , we obtain a geodesic tree
connecting at least $r+1$ singular points of $\Sig$.\\

\noindent Let $C_3$ denote the new tree. For any point $p$ of $C_3$, the length of the unique path from $p_1$ to $p$ along $C_3$ is the distance
$\di(p_1,p)$. This property allows us to conclude by induction.  \carre

\subsection{Proof of Proposition \ref{SSprA}}

Since the complement of a tree in the sphere $\S^2$ is a topological disk, on a spherical flat surface,  any geodesic tree whose vertex set is
the set of singular points is automatically an erasing tree.  Therefore, the arguments for flat surfaces with erasing forest can be applied in
the case of spherical flat surfaces. \\

\noindent Let $\TRSn$ denote the set of triangulations of $\S^2$ whose vertex set is $\{p_1,\dots,p_n\}$ modulo isotopy relative to
$\{p_1,\dots,p_n\}$. Given a triangulation $\Tt$ of $\S^2$, which represents an equivalence class in $\TRSn$,  let $\U_\Tt$ denote the subset of
$\TSS$ consisting of pairs $([(\Sig,\phi)],e^{\imath\theta})$, such that $\phi(\Tt)$ is a geodesic triangulation of $\Sig$.\\

\noindent Pick a tree $\A$ in the $1$-skeleton of  $\Tt$ whose vertex set is $\{p_1,\dots,p_n\}$, for any $([(\Sig,\phi)],e^{\imath\theta})$ in
$\U_\Tt$, $\phi(\A)$ is a geodesic erasing tree of $\Sig$. Therefore, we can identify $\U_\Tt$ to an open subset in
$\mathcal{T}^{\mathrm{et}}(\S^2,\A)$. From the proof of Theorem \ref{ETthB}, we get a map $\DS{\Psi_{\Tt,\A} :\U_\Tt\lra \C^{4n-7}}$ which is
injective, such that $\Psi_{\Tt,\A}(\U_\Tt)$ is an open subset of the solution space $\mathrm{V}_{\Tt,\A}$ of a system of linear equations
$\mathbf{S}_{\Tt,\A}$. The maps $\Psi_{\Tt,\A}$ are local charts of $\TSS$. Since in this case $\dim_\C\mathrm{V}_{\Tt,\A}=n-2$, it follows that
$\dim_\C \TSS=n-2$. It is worth noticing that $\Psi_{\Tt,\A}$ is only defined up to a rotation.\\

Given two triangulations $\Tt_1,\Tt_2$ representing distinct equivalence classes in $\TRSn$, let $([(\Sig,\phi)],e^{\imath\theta})$ be a point
in $\U_{\Tt_1}\cap\U_{\Tt_2}$, and let $\T_1$, $\T_2$ be the geodesic triangulations of $\Sig$ corresponding to $\Tt_1$, and $\Tt_2$
respectively. Choose a tree $\A_1$ (resp. $\A_2$) in $\Tt_1$ (resp. $\Tt_2$) which connects all the points in $\{p_1,\dots,p_n\}$, and let
$\Psi_{\Tt_1,\A_1}$ and $\Psi_{\Tt_2,\A_2}$ be the two local charts of $\TSS$ corresponding.\\

\noindent Let $e$ be an edge of $\T_2$ which is not contained in $\T_1$, and  $\P_e$ be the developing polygon of $e$ with respect to $\T_1$
(see Lemma \ref{lmA10}). By construction, the complex number associated to the edge $e$ in the local chart $\Psi_{\Tt_2,\A_2}$ can be written as
a linear function of complex numbers associated to edges of $\T_1$, which correspond the segments in the boundary of $\P_e$ in the local chart
$\Psi_{\Tt_1,\A_1}$. Since the roles of $\T_1$ and $\T_2$ in this reasoning can be interchanged, we deduce that the coordinate change between
$\Psi_{\Tt_1,\A_1}$ and $\Psi_{\Tt_2,\A_2}$ is a linear isomorphism of $\C^{4n-7}$ which sends $\mathrm{V}_{\Tt_1,\A_1}$ onto
$\mathrm{V}_{\Tt_2,\A_2}$. Thus, we can conclude that $\TSS$ is a flat complex affine manifold of dimension $n-2$.\\

We can define a map $\DS{\Xi:\TSS \ra \Tsn}$ by associating to each point $([(\Sig,\phi)],e^{\imath\theta})\in \TSS$, the equivalence class of
$(\Sig,\phi)$ in $\Tsn$, where $\Sig$ is now considered as a Riemann surface. Using the same arguments as in Proposition \ref{ETprB3}, we can
see that $\Xi$ is continuous. Clearly, the map $\Xi$ is equivariant with respect to the actions of $\Gsn$ on $\TSS$, and $\Tsn$, it follows that
the action of $\Gsn$ on $\TSS$ is properly discontinuous. The proof of Proposition \ref{SSprA} is now complete. \carre

\rem The map $\Xi$ can be defined on $\TSSv$. Using Proposition 2 \cite{Troy}, we can show that the restriction of $\Xi$ onto $\TSSvi$, the
subspace of $\TSSv$ consisting of surfaces of unit area, is a bijection. Therefore, we can consider $\TSS$ as a $\C^*$-bundle over $\Tsn$.

\section{Volume form on moduli space of flat surfaces of genus zero}\label{SSVolSec}

In this section, we will give the proof of Theorem \ref{SSthB}. Set $N_1=4n-7$, $N_2=3n-5$, and let $\Tt$ be a triangulation of $\S^2$
representing an equivalence class in $\TRSn$. Let $\A$ be a tree contained in $\Tt$, which connects all the points in $\{p_1,\dots,p_n\}$. Let
$\Psi_{\Tt,\A}$ be the local chart associated to $(\Tt,\A)$, which is defined on the set $\U_\Tt$. Following the method in Section
\ref{ETVolSec}, we then get a volume form $\mu_{\Tt,\A}$ on $\U_\Tt$. Recall that there is a normalized linear map $\displaystyle{\ATtA:
\C^{N_1} \ra \C^{N_2}}$ associated to the local chart $\Psi_{\Tt,\A}$, and by definition,
$\displaystyle{\mu_{\Tt,\A}=\Psi_{\Tt,\A}^*\nu_{\Tt,\A}}$, where $\nu_{\Tt,\A}$ is the volume form on $\ker \ATtA$ which is induced by the
Lebesgue measures on $\C^{N_1}$ and $\C^{N_2}$ via the following exact sequence

\begin{equation}\label{ExSeq1b}
0\lra \ker \ATtA \hookrightarrow \C^{N_1}\overset{\ATtA}{\lra}\C^{N_2}\lra 0
\end{equation}

\noindent  The following proposition shows that the volume form $\mu_{\Tt,\A}$ does not depend on the choice of $\A$.\\

\begin{proposition}\label{SSPrB1}

Let $\Tt$ be a triangulation representing an equivalence class in $\TRSn$. Let $\A_1,\A_2$ be two trees contained in the $1$-skeleton of $\Tt$,
each of which connects all the points in $\{p_1,\dots,p_n\}$. Let $\mathbf{A}_{\Tt,\A_1}$ and $\mathbf{A}_{\Tt,\A_2}$ denote the linear maps
from $\C^{N_1}$ onto $\C^{N_2}$ corresponding to $\A_1$, and $\A_2$ respectively. Let $\nu_i,\; i=1,2$ denote the volume form on
$\ker\mathbf{A}_{\Tt,\A_i}$ which is defined by the exact sequence (\ref{ExSeq1b}). Let $\mathbf{H}=\Psi_{\Tt,\A_2}\circ \Psi_{\Tt,\A_1}^{-1}$
be the coordinate change between $\Psi_{\Tt,\A_1}$, and $\Psi_{\Tt,\A_2}$, then we have $\DS{\mathbf{H}^*\nu_2=\nu_1}$
\end{proposition}

\noindent To show that the volume form $\mu_{\Tt,\A}$ actually does not depend on the choice of $\Tt$, we need the following

\begin{theorem}\label{SSThD}

Let $\Sig$ be a spherical flat surface. If $\T_1$ and $\T_2$ are two geodesic triangulations of $\Sig$ whose sets of vertices coincide, and
contain the set of singularities of $\Sig$, then $\T_1$ and $\T_2$ are related (\ie one can be transformed into the other by elementary
moves).\\

\end{theorem}

Theorem \ref{SSthB} follows directly from Proposition \ref{SSPrB1}, and Theorem \ref{SSThD}, since by Lemma \ref{ETlmB7} we know that, if $\T_2$
is obtained from $\T_1$ by an elementary move then the volume forms corresponding to $\T_1$ and $\T_2$ coincide. Theorem \ref{SSThD} is of
course a consequence of the fact that any geodesic triangulation of a spherical flat surface whose vertex set coincides with the set of
singularities can be transformed into a Delaunay triangulation by elementary moves. In Appendices, Section \ref{prfSSThD}, we give a proof of
Theorem \ref{SSThD} using similar ideas to the proof of Theorem \ref{ETthF}. The remainder of this section is devoted to the proof of
Proposition \ref{SSPrB1}.\\

\subsection{Cutting and gluing}

Let us consider pairs $(\Sig_0,\T_0)$ where
\begin{itemize}

\item[-] $\Sig_0$ is a flat surface homeomorphic to a closed disk, with geodesic boundary, and having no singularities in the interior.

\item[-] $\T_0$ is a triangulation of $\Sig_0$ by geodesic segments whose vertex set is contained in the boundary of $\Sig^0$.

\item[-] The edges of $\T_0$ on the boundary of $\Sig_0$ are paired up. Two edges in a pair have the same length.

\end{itemize}

\noindent We will call such a pair a {\em well triangulated flat disk}. Remark that a geodesic tree $A$ contained in $\T$ which connects all the
singular points of $\Sig$ gives rise to a well triangulated flat disk, which is obtained by slitting open $\Sig$ along $A$.\\

\noindent Let $\Tt,\A_1,\A_2$  be as in Proposition \ref{SSPrB1}. Given a point  $([(\Sig,\phi)], e^{\imath\theta})$ in $\U_\Tt$, let $\T$ be
the geodesic triangulation of $\Sig$ corresponding to $\Tt$, and $A_1,A_2$ be the geodesic trees corresponding to $\A_1,\A_2$ respectively. Let
$\Sig_0^1$ and $\Sig_0^2$ denote the flat surface with geodesic boundary obtained by slitting open the surface $\Sig$ along the trees $A_1$ and
$A_2$ respectively. Let $\T_0^1$ (resp. $\T_0^2$) denote the geodesic triangulation of $\Sig_0^1$ (resp. $\Sig_0^2$) which is induced by $\T$.
By definition, $(\Sig_0^1,\T_0^1)$ and $(\Sig_0^2,\T_0^2)$ are well triangulated flat disks. Consider the following the following operation:

\begin{itemize}

\item[$\bullet$] Choose a pair of edges $(a,\bar{a})$ of $\T_0$ in the boundary of $\Sig_0$, and an edge $b$ in the interior of $\Sig_0$ so that
$a$ and $\bar{a}$ do not belong to the same connected component of $\Sig_0\setm b$.

\item[$\bullet$] Cut $\Sig_0$ along $b$, then glue two the sub-disks  by identifying $a$ to $\bar{a}$.

\end{itemize}

\noindent Clearly, by this operation, we get another pair $(\Sig'_0,\T'_0)$ with is also a well triangulated flat disk. We will call this
operation the {\em cutting-gluing operation}. We have:

\begin{lemma}\label{SSLmB1}
The pair $(\Sig_0^2,\T_0^2)$ can be obtained from $(\Sig^1_0,\T^1_0)$ by a sequence of cutting-gluing operations.
\end{lemma}

\dem  First, let us prove that $A_2$ is obtained from $A_1$ by replacing edges of $A_1$ by edges of $A_2$ one by one, successively. Indeed, let
$e$ be an edge of $\T$ which is contained in $A_2$, but not in $A_1$. Let $v_1$ and $v_2$ denote the endpoints of $e$, by assumption, there is a
path $c$ in $A_1$ which joins $v_1$ to $v_2$. The union of $c$ and $e$ is then a cycle in $\T$, therefore there exists an edge $e'$ in $c$,
different from $e$, which does not belong to $A_2$. Replacing $e'$ by $e$, we get a new tree which connects all the singular points, and
contains one more common edge with $A_2$ than $A_1$. The claim follows by induction.\\

\noindent  Now, we just need to observe that the operation of replacing $e'$ by $ e$ corresponds to a cutting-gluing operation on the well
triangulated flat disk corresponding to the tree $A_1$, and the lemma follows. \carre


\subsection{Increased exact sequence}

Given a well triangulated flat disk $(\Sig_0,\T_0)$ arising from a tree $A$ contained in the $1$-skeleton of $\T$ which connects all the
singular points of $\Sig$, we have a system of linear equations $\mathbf{S}_0$ defined as in \ref{ETlocch}, and a normalized complex linear map
$\mathbf{A}_0: \C^{N_1} \lra \C^{N_2}$ corresponding. Let $a_1,\dots,a_{N_2}$ denote the row vectors of $\mathbf{A}_0$. We know that
$\mathrm{rk}(\mathbf{A}_0)=N_2$. Now, choose an edge $e_0$ of $\T_0$ which is contained inside $\Sig_0$, and cut $\Sig_0$ along $e_0$, we then
get two flat disks $\mathbf{D}_1,\mathbf{D}_2$ with the geodesic triangulations $\T_1, \T_2$ respectively. Again, we can associate to each edge
$e$ of $\T_1$ and $\T_2$ a complex number $z(e)$ to get a vector $\hat{Z}_0$ in $\C^{N_1+1}$. We also have a system of $N_2$ linear equations
arising from the equations of $\mathbf{S}_0$. We add to this system the following equation

$$z(e'_0)+z(e''_0)=0$$

\noindent where $e'_0$ and $e''_0$ are the edges of $\T_1$ and $\T_2$ which arise from $e_0$ (with appropriate orientation). Let
$\hat{\mathbf{S}}_0$  denote the new system, and $\hat{\mathbf{A}}_0: \C^{N_1+1} \lra \C^{N_2+1}$ denote the normalized linear map
corresponding. We will call $\hat{\mathbf{A}}_0$ the {\em increased linear map} of $\mathbf{A}_0$ associated to the splitting along $e_0$.\\

\noindent By construction, there exists injective linear maps $\mathbf{J}_1: \C^{N_1}\lra \C^{N_1+1}$, $\mathbf{J}_2:\C^{N_2}\lra \C^{N_2+1}$,
and an isomorphism $\mathbf{I}:\ker \mathbf{A}_0 \lra \hat{\mathbf{A}}_0$ so that the following diagram is commutative:

\begin{equation}\label{Diag1}
\begin{array}{ccccccccc}
  0 & \lra & \ker\mathbf{A}_0 & \overset{\iota}{\hookrightarrow} & \C^{N_1} & \overset{\mathbf{A}_0}{\lra} & \C^{N_2} & \lra & 0 \\
   &  & \downarrow \mathbf{I} &  & \downarrow \mathbf{J}_1 &  & \downarrow \mathbf{J}_2 &  &  \\
  0&\lra & \ker\hat{\mathbf{A}}_0& \overset{\hat{\iota}}{\hookrightarrow} & \C^{N_1+1} & \overset{\hat{\mathbf{A}}_0}{\lra}& \C^{N_2+1}&\lra& 0 \\
\end{array}
\end{equation}

\noindent As usual, let $\lambda_{2k}$ denote the Lebesgue measure of $\C^{k}, \; \forall k\in \N$. Let $\nu$ denote the volume form on
$\ker\mathbf{A}_0$ which is induced by $\lambda_{2N_1}$ and $\lambda_{2N_2}$ via the upper exact sequence in (\ref{Diag1}), and $\hat{\nu}$ is
the volume form on $\ker\hat{\mathbf{A}}_0$ which is induced by $\lambda_{2(N_1+1})$ and $\lambda_{2(N_2+1)}$ via the lower exact sequence in
(\ref{Diag1}). We can choose a numbering of the edges of $\T_0,\T_1,\T_2$ so that

$$\mathbf{J}_1((z_1,\dots,z_{N_1}))=(z_1,\dots,z_{N_1},-z_1), \text{ and } \mathbf{J}_2((z_1,\dots,z_{N_2}))=(z_1,\dots,z_{N_2},0).$$

\noindent Let $a_i$ (resp. $\hat{a}_i$), $i=1,\dots,N_2$ denote the $i$-th row of $\mathbf{A}_0$ (resp. $\hat{\mathbf{A}}_0$), we consider $a_i$
(resp. $\hat{a}_i$) as complex linear map from $\C^{N_1}$ (resp. $\C^{N_1+1}$) into $\C$. By assumption, we have $a_i=\mathbf{J}_1^*\hat{a}_i$.
Let $\mathbf{h}: \C^{N_1+1}\ra \C$ be the complex linear map corresponding to the last row of $\hat{\mathbf{A}}_0$, by assumption,  we have
$\DS{\mathbf{h}((z_1,\dots,z_{N_1+1}))=z_1+z_{N_1}}$, and the following sequence

\begin{equation}\label{ExSeq1c}
0\lra \C^{N_1}\overset{\mathbf{J}_1}{\hookrightarrow} \C^{N_1+1}\overset{\mathbf{h}}{\lra} \C \lra 0
\end{equation}

\noindent is exact. Let $\lambda'_{2N_1}$ be the volume form induced by $\lambda_{2(N_1+1)}$ and $\lambda_2$ via the exact sequence
(\ref{ExSeq1c}). Let $\nu'$ be the volume form on $\ker\mathbf{A}_0$ which is induced by $\lambda'_{2N_1}$ and $\lambda_{2N_2}$ via the upper
exact sequence in (\ref{Diag1}). We have

\begin{lemma}\label{SSLmB2}
$\mathbf{I}^*\hat{\nu}$ is equal to $\nu'$ .
\end{lemma}

\dem Let $\hat{\eta}$ be a $2(N_1-N_2)$-real form on $\C^{N_1+1}$ such that

$$\hat{\eta}\wedge(\Re(\hat{a}_1)\wedge\Im(\hat{a}_1))\wedge\dots\wedge(\Re(\hat{a}_{N_2}) \wedge \Im(\hat{a}_{N_2}))\wedge(\Re(\mathbf{h}) \wedge \Im
(\mathbf{h}))= \lambda_{2(N_1+1)}.$$

\noindent By definition, we have $\hat{\nu}=\hat{\iota}^*\hat{\eta}$, where $\hat{\iota}: \ker \hat{\mathbf{A}}_0 \hookrightarrow \C^{N_1+1}$ is
the natural embedding. On the other hand, by definition, we have

$$\lambda'_{2N_1}=\mathbf{J}_1^*(\hat{\eta}\wedge(\Re(\hat{a}_1)\wedge\Im(\hat{a}_1)) \wedge\dots\wedge(\Re(\hat{a}_{N_2}) \wedge
\Im(\hat{a}_{N_2})).$$

\noindent Therefore

$$\lambda'_{2N_1}=\mathbf{J}_1^*\hat{\eta}\wedge (\Re(a_1)\wedge\Im(a_1)) \wedge\dots\wedge(\Re(a_{N_2}) \wedge \Im(a_{N_2})),$$

\noindent and by definition, the volume form $\nu'$ is equal to $\iota^*(\mathbf{J}_1^*\hat{\eta})$. Since $\mathbf{J}_1\circ \iota=
\hat{\iota}\circ\mathbf{I}$, it follows $\nu'=\mathbf{I}^* (\hat{\iota}^*\hat{\eta}) =\mathbf{I}^*\hat{\nu}$.\carre

\begin{corollary}\label{SScorB2}
We have $\DS{\mathbf{I}^*\hat{\nu}=c_0\nu}$, where $c_0$ is a constant which does not depend on the choice of the edge $e_0$.
\end{corollary}

\dem Set $\displaystyle{c_0=\frac{\lambda'_{2N_1}}{\lambda_{2N_1}}}$, then we have $\nu'=c_0\nu$. By lemma \ref{SSLmB2}, we deduce that
$\mathbf{I}^*\hat{\nu}=c_0\nu$. Now, the splitting of $\Sig_0$ along another edge of $\T_0$ inside $\Sig_0$ corresponds to a permutation of
coordinates in $\C^{N_1+1}$, hence the exact sequence (\ref{ExSeq1c}) gives the same volume form $\lambda'_{N_1}$ on $\C^{N_1}$, and the lemma
follows.\carre


\subsection{Proof of Proposition \ref{SSPrB1}}

By Lemma \ref{SSLmB1}, it suffices to consider the case where $(\Sig^2_0,\T^2_0)$ is obtained from $(\Sig^1_0,\T^1_0)$ by only one
cutting-gluing operation. We can then assume $(\Sig_0^2,\T_0^2)$ is obtained by cutting $\Sig^1_0$ along an edge $e_1$ of $\T_0^1$ inside
$\Sig_0^1$, and gluing a pair $(e,\bar{e})$ of edges of $\T_0^1$ in $\partial\Sig_0^1$ which gives raise to an edge $e_2$ of $\T_0^2$.\\

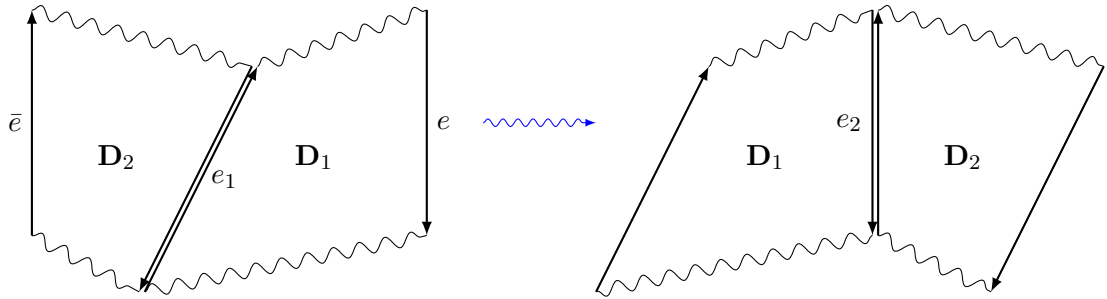
\begin{center}
\begin{tikzpicture}[>=latex, scale=0.75]
\draw[->, thick] (0,-2) -- (0,2); \draw[decorate, decoration=snake] (0,2) -- (3.9,1); \draw[decorate, decoration=snake] (4,1) -- (7,2);
\draw[decorate, decoration=snake] (0,-2) -- (1.9,-3); \draw[decorate, decoration=snake] (2,-3) -- (7,-2); \draw[->,thick] (3.9,1)-- (1.9,-3);
\draw[->,thick] (2,-3) -- (4,1); \draw[->,thick] (7,2) -- (7,-2);

\draw[<-, thick] (12,1) -- (10,-3); \draw[decorate, decoration=snake] (12,1) -- (14.9,2); \draw[decorate, decoration=snake] (15,2) -- (19,1);
\draw[->, thick] (15,-2) -- (15,2); \draw[->, thick] (14.9,2) -- (14.9, -2); \draw[decorate, decoration=snake] (14.9,-2) -- (10,-3);
\draw[decorate, decoration=snake] (15,-2) -- (17,-3); \draw[->, thick] (19,1) -- (17,-3);

\draw (0,0) node[left] {$\bar{e}$}; \draw (1.5,-1) node[above] {$\mathbf{D}_2$}; \draw (3,-1) node[right] {$e_1$}; \draw (5,-1) node[above]
{$\mathbf{D}_1$}; \draw (7,0) node[right] {$e$}; 

\draw (13,-1) node[above] {$\mathbf{D}_1$}; \draw (16.5,-1) node[above] {$\mathbf{D}_2$}; \draw (14.9,0) node[left] {$e_2$};


\draw[blue, ->, decorate, decoration={snake, amplitude=0.5mm, segment length=2mm, post length=1mm}] (8,0) -- (10,0);
\end{tikzpicture}

\end{center}

\noindent Let $\hat{\mathbf{A}}_1,\hat{\mathbf{A}}_2: \C^{N_1+1}\lra \C^{N_2+1}$ denote the linear maps corresponding to the splitting of
$\Sig_0^1$ and $\Sig_0^2$ along $e_1$ and $e_2$ respectively. There exist the isomorphisms $\hat{\mathbf{F}}:\C^{N_1+1}\rightarrow \C^{N_1+1},
\hat{\mathbf{G}}:\C^{N_2+1}\rightarrow\C^{N_2+1}$ such that the following diagram is commutative

\begin{equation}\label{Diag2}
\begin{array}{ccccccccc}
  0 & \lra & \ker \hat{\mathbf{A}}_1 & \overset{\hat{\iota}_1}{\lra} & \C^{N_1+1} & \overset{\hat{\mathbf{A}}_1}{\lra} & \C^{N_2+1} & \lra & 0 \\
   &  & \downarrow \hat{\mathbf{H}} &  & \downarrow \hat{\mathbf{F}} &  & \downarrow \hat{\mathbf{G}} &  &  \\
  0 & \lra & \ker \hat{\mathbf{A}}_2&\overset{\hat{\iota}_2}{\lra} & \C^{N_1+1} &\overset{\hat{\mathbf{A}}_2}{\lra} & \C^{N_2+1} & \lra & 0\\
\end{array}
\end{equation}

\noindent where $\hat{\mathbf{H}}$ is the isomorphism induced by $\hat{\mathbf{F}}$. Let $k$ be the number of edges of $\T_1$ which are
contained in one of the two disks obtained from the splitting of $\Sig_0^1$ along $e_1$, and $\theta$ be the angle of the rotation
$\orth(\gamma_{(e,\bar{e})})$, where $\gamma_{(e,\bar{e})}$ is a closed curve on $\Sig$ corresponding to a curve in $\Sig_0^1$ joining the
midpoint of $e$ to the midpoint of $\bar{e}$. With appropriate numberings of the edges of $\T_0^1$, and the edges of $\T_0^2$, we have

$$\hat{\mathbf{F}}=\left(%
\begin{array}{cc}
  e^{\imath\theta}\Id_k & 0 \\
  0 & \Id_{N_1+1-k} \\
\end{array}%
\right).$$

\noindent Consequently, $\hat{\mathbf{G}}$ is a diagonal matrix in $\mathbf{M}_{N_2+1}(\C)$ whose diagonal entries are either $1$ or
$e^{\imath\theta}$. Clearly, we have

$$|\det\hat{\mathbf{F}}|=|\det \hat{\mathbf{G}}|=1.$$

\noindent It follows that  $\hat{\mathbf{H}}^*\hat{\nu}_2=\hat{\nu}_1$, where $\hat{\nu}_i$ is the volume form on $\ker\hat{\mathbf{A}}_i$ which
is induced by the Lebesgue measures on $\C^{N_1+1}$, and $\C^{N_2+1}$.\\

Let $\mathbf{I}_i : \ker\mathbf{A}_i\ra \ker \hat{\mathbf{A}}_i, \; i=1,2,$ denote the isomorphism in (\ref{Diag1}) corresponding to the
splitting of $\Sig_i$ along $e_i$. We have the following commutative diagram

\begin{equation}\label{Diag3}
\begin{array}{ccc}
  \ker\mathbf{A}_1 & \overset{\mathbf{I}_1}{\lra} & \ker\hat{\mathbf{A}}_1 \\
  \downarrow \mathbf{H}&  & \downarrow \hat{\mathbf{H}} \\
  \ker \mathbf{A}_2 & \overset{\mathbf{I}_2}{\lra} & \ker\hat{\mathbf{A}}_2 \\
\end{array}
\end{equation}

\noindent By Corollary \ref{SScorB2}, we know that $\DS{\frac{{\mathbf{I}_1}^* \hat{\nu}_1}{\nu_1}=\frac{{\mathbf{I}_2}^* \hat{\nu}_2}{\nu_2}}$,
and the proposition follows.\carre

\section{Comparison with complex hyperbolic volume form}\label{SSCompSec}

\subsection{Definitions}

In this section, we assume that all the angles $\alpha_1,\dots,\alpha_n$ are less than $2\pi$. Put $\kappa_i=2\pi-\alpha_i, \; i=1,\dots,n$, we
have $\displaystyle{\kappa_1+\dots+\kappa_n=4\pi}$. Following Thurston \cite{Thu}, we denote by $C(\kappa_1,\dots,\kappa_n)$ the moduli space of
spherical flat surface having $n$ singularities with cone angles $\alpha_1,\dots,\alpha_n$, or equivalently, with curvatures $\kappa_1,\dots,
\kappa_n$, up to homothety. In \cite{Thu}, Thurston proves that $C(\kappa_1,\dots,\kappa_n)$ admits a complex hyperbolic metric structure with
finite volume, and the metric closure of $C(\kappa_1,\dots,\kappa_n)$ has cone manifold structure.\\

\noindent The complex hyperbolic metric provides a volume form $\mu_\mathrm{Hyp}$ on $C(\kappa_1,\dots,\kappa_n)$. On the other hand, the volume
form $\mu_\mathrm{Tr}$ gives another volume form $\hat{\mu}^1_\mathrm{Tr}$ on $C(\kappa_1,\dots,\kappa_n)$ which is defined as follows:

\begin{itemize}

\item[-] First, we identify $C(\kappa_1,\dots,\kappa_n)$ to the subset $\MSSvi$ of all surfaces of area $1$ in $\MSSv$. Let $f:\MSS \ra \R^+$ be
the function which associates to a pair $(\Sig,\theta)$ in $\MSS$ the area of $\Sig$. The space $\MSSvi$ can be considered as the quotient of
the locus $f^{-1}(\{1\})$ by the action of $\S^1$.

\item[-] By Proposition \ref{SSprA}, we know that $\MSS$ is a complex orbifold, let $\J$ denote the complex structure of $\MSS$. Let
$\displaystyle{\rho : f^{-1}(\{1\}) \ra f^{-1}(\{1\})/\S^1=\MSSvi}$ denote the natural projection. We define the volume form
$\hat{\mu}^1_\mathrm{Tr}$ on $\MSSvi$ to be the one such that:

\begin{equation*}
\rho^*\hat{\mu}^1_\mathrm{Tr}\wedge df\wedge (df\circ \J)=\mu_\mathrm{Tr}
\end{equation*}

\end{itemize}

\subsection{Local formulae}

First, we recall the construction of local charts for $C(\kappa_1,\dots,\kappa_n)$ as presented in \cite{Thu}, and  consequently the definition
of $\mu_\mathrm{Hyp}$. Given a surface $\Sig$ in $\MSSvi$, we consider $\Sig$ as a point in $C(\kappa_1,\dots,\kappa_n)$. Let $\T$ be a
triangulation of $\Sig$ by geodesic segments whose set of vertices is the set of singular points. Choose a singular point of $\Sig$ and denote
this point by $x_{last}$. We will call all the edges of $\T$ which contain $x_{last}$ as an endpoint {\em followers }. Pick a tree $\tilde{A}$
in $\T$ which connects all other singular points of $\Sig$, and call the edges of this tree {\em leaders}. The remaining edges of $\T$ are also
called {\em followers}.\\

\noindent Using a developing map, one can associate to each of the leaders  a complex number, there are $n-2$ of them. Let $(z_1,\dots,z_{n-2})$
denote those complex numbers. The same developing map also defines an associated complex number for each of the followers, but these numbers can
be calculated from those associated to leaders by complex linear functions. Thus, the complex numbers associated to leaders determine a local
coordinate system $\varphi : \mathrm{U} \lra \MSS$ for $\MSS$ in a neighborhood of  $(\Sig,1)$, where $\mathrm{U}$ is a neighborhood of
$(z_1,\dots,z_{n-2})$ in $\C^{n-2}$. Consequently, a neighborhood of $\Sig$ in $C(\kappa_1,\dots,\kappa_n)$ is then identified to an open set of
$\mathbb{P}\C^{n-3}$ which contains $[z_1:\dots:z_{n-2}]$.\\

\noindent If we add to $\tilde{A}$ a follower which joins $x_{last}$ to the tree $\tilde{A}$, then we have an erasing tree $A$ on $\Sig$. We can
then construct a local chart $\Psi_{\Tt,\A}$ for $\MSS$ from $\T$ and $A$ associated to a complex linear surjective map $\ATt: \C^{N_1}\ra
\C^{N_2}$ is determined by the tree $A$. In this local chart,  $\mu_\mathrm{Tr}$ is identified to the volume form on $\ker \ATt$ which is
induced by the exact sequence (\ref{ExSeq1b}).\\

\noindent Now, observe that the map $\Psi_{\Tt,\A}\circ \varphi$ is the restriction to an open subset of $\C^{n-2}$ of an isomorphism from
$\C^{n-2}$ to $\ker \ATt$. Therefore, the following sequence is exact

$$0\lra \C^{n-2} \overset{\Psi_{\Tt,\A}\circ\varphi}{\lra} \C^{N_1}\overset{\ATt}{\lra} \C^{N_2}\lra 0.$$

\noindent Consequently, in the local chart $\varphi$, we have

$$\mu_\mathrm{Tr}=c\lambda_{2(n-2)},$$

\noindent where $\lambda_{2(n-2)}$ is the Lebesgue measure of $\C^{n-2}$, and $c$ is a constant.\\

\noindent In the local chart $\varphi$, the area function $f$ on $\MSS$ is expressed as a Hermitian form $\Her$. More precisely, if $v\in
\C^{n-2}$ is a vector such that $\displaystyle{\varphi(v)=(\Sig,\theta)\in \MSS}$ then $\displaystyle{ f((\Sig,\theta))=\mathbf{Area}(\Sig)=
\transpose{\overline{v}}\Her v}$. It is proven in \cite{Thu} that $\Her$ is of signature $(1,n-3)$. Changing the basis and the sign of $\Her$,
we can assume that

$$ \Her =\left(%
\begin{array}{cc}
  \Id_{n-3} &  0 \\
   0 & -1 \\
\end{array}%
\right)$$

\noindent Thus we can write $\displaystyle{f (z_1,\dots,z_{n-2})=|z_1|^2+\dots+|z_{n-3}|^2-|z_{n-2}|^2}$. Note that by these changes, the
vectors of $\C^{n-2}$ representing surfaces in $\MSSvi$ are contained in the set $\mathbf{Q}_1=f^{-1}(\{-1\})$, and we still have
$\mu_\mathrm{Tr}=c_0\lambda_{2(n-2)}$ with $c_0$ a constant.\\

We use the symbol $\langle,\rangle$ to denote the scalar product defined by Hermitian form $\Her$. By definition $f(Z)=\langle Z,Z \rangle,\;
\forall Z \in \C^{n-2}$. Let $\J$ be the natural complex structure of $\C^{n-2}$, that is $\J(z_1,\dots,z_{n-2})=(\imath z_1,\dots,\imath
z_{n-2})$, and  $\eta$ be the real symmetric form induced by $\langle,\rangle$, that is

$$\eta(X,Y)=\mathrm{Re}\langle X, Y \rangle.$$

\noindent Let $Z$ be a vector in $\mathbf{Q}_1$ which represents a surface in $\MSSvi$. The tangent space of $\mathbf{Q}_1/\S^1$ at the orbit
$\S^1\cdot Z $ is naturally identified to the orthogonal complement of $Z$ with respect to $\langle,\rangle$, we denote this space by $Z^\bot$.
The restriction of $\langle,\rangle$ on $Z^\bot$ is a definite positive Hermitian form, which determines the complex hyperbolic metric on
$\MSSvi=C(\kappa_1,\dots,\kappa_n)$. We have

$$df=(\bar{z_1}dz_1+\dots+\bar{z}_{n-3}dz_{n-3}-\bar{z}_{n-2}dz_{n-2})+(z_1d\bar{z}_1+\dots+z_{n-3}d\bar{z}_{n-3}-z_{n-2}d\bar{z}_{n-2}),$$

\noindent and

$$df\circ \J=
\imath(\bar{z_1}dz_1+\dots+\bar{z}_{n-3}dz_{n-3}-\bar{z}_{n-2}dz_{n-2})-\imath(z_1d\bar{z}_1+\dots+z_{n-3}d\bar{z}_{n-3}-z_{n-2}d\bar{z}_{n-2}).$$

\noindent Note that both $df$ and $df\circ \J$ are invariant by the action of $\S^1$. Put

$$U_k=(0,\dots,0,\overset{(k)}{\bar{z}_{n-2}},0,\dots,\bar{z}_k), \text{ and} \; V_k=\J\cdot U_k=\imath U_k, \text{ for} \; k=1,\dots,n-3.$$

\noindent One can easily check that $\{U_1,V_1,\dots,U_{n-3},V_{n-3}\}$ span $Z^\bot$ as a real vector space. We consider
$\{U_1,V_1,\dots,U_{n-3},V_{n-3}\}$ as a basis of the tangent space of $\MSSvi$ at $\varphi(Z)$. We know that the restriction of the symmetric
form $\eta$ on $Z^\bot$ defines a Riemannian metric. Let $U^*_k,V^*_k$ denote the $\R$-linear $1$-forms dual to $U_k$ and $V_k$ respectively
with respect to $\eta$. We have

$$U^*_k=\frac{1}{2}[(z_{n-2}dz_k-z_kdz_{n-2})+(\bar{z}_{n-2}d\bar{z}_{k}-\bar{z}_kd\bar{z}_{n-2})],$$

\noindent and

$$V^*_k=\frac{-\imath}{2}[(z_{n-2}dz_k-z_kdz_{n-2})-(\bar{z}_{n-2}d\bar{z}_{k}-\bar{z}_kd\bar{z}_{n-2})].$$

\noindent We can consider $\{U^*_1,V^*_1,\dots,U^*_{n-3},V^*_{n-3}\}$ as a basis of the cotangent space of $\MSSvi$ at $\varphi(Z)$. Let $\rho$
be the projection from $\mathbf{Q}_1$ to $\mathbf{Q}_1/\S^1$. By definition, the volume form  $\hat{\mu}^1_\mathrm{Tr}$ on $\mathbf{Q}_1/\S^1$
verifies

\begin{equation}\label{VolFSSvieq}
\rho^*\hat{\mu}^1_\mathrm{Tr}\wedge df\wedge (df\circ \J)=d\lambda_{2(n-2)}=(\frac{\imath}{2})^{n-2}dz_1d\bar{z}_1\dots dz_{n-2}d\bar{z}_{n-2}
\end{equation}

\noindent Since $df$ and $df\circ \J$ are invariant by the action of $\S^1$, the volume form $\hat{\mu}^1_\mathrm{Tr}$ is well defined by this
condition. We will express $\hat{\mu}^1_\mathrm{Tr}(\S^1\cdot Z)$ in terms of $U^*_k,V^*_k,\; k=1,\dots,n-3$.

\begin{lemma}\label{SSlmC1}
We have

$$\hat{\mu}^1_\mathrm{Tr}(\S^1\cdot Z)=\frac{c_0}{|z_{n-2}|^{2(n-4)}}(U^*_1\wedge V^*_1)\wedge\dots\wedge(U^*_{n-3}\wedge V^*_{n-3}),$$

\noindent where $\DS{c_0=\frac{\mu_\mathrm{Tr}}{\lambda_{2(n-2)}}}$.

\end{lemma}

\dem Set

\begin{itemize}

\item[$\bullet$] $X_k=z_{n-2}dz_k-z_kdz_{n-2}, \; \overline{X}_k=\bar{z}_{n-2}d\bar{z}_k-\bar{z}_kd\bar{z}_{n-2}, \; k=1,\dots,n-3$, and

\item[$\bullet$] $ X=\bar{z}_1dz_1+\dots+\bar{z}_{n-3}dz_{n-3}-\bar{z}_{n-2}dz_{n-2}, \; \overline{X}=
z_1d\bar{z}_1+\dots+z_{n-3}d\bar{z}_{n-3}-z_{n-2}d\bar{z}_{n-2}.$

\end{itemize}

\noindent We have

\begin{itemize}

\item[$\bullet$] $\displaystyle{ U^*_k\wedge V^*_k =\frac{-\imath}{4}(X_k+\overline{X}_k)\wedge
(X_k-\overline{X}_k)=\frac{\imath}{2}X_k\wedge\overline{X}_k}, \; k=1,\dots,n-3$, and

\item[$\bullet$] $\displaystyle{df\wedge (df\circ\J)=2\imath X\wedge\overline{X}}$.

\end{itemize}

\noindent Thus

$$\begin{array}{ll}
 & (U^*_1\wedge V^*_1\wedge\dots\wedge U^*_{n-3}\wedge V^*_{n-3})\wedge df \wedge (df\circ\J)\\
 =& -(\frac{\imath}{2})^{n-4}X_1\wedge \overline{X}_1\wedge\dots\wedge X_{n-3}\wedge\overline{X}_{n-3}\wedge X\wedge\overline{X}\\
=& -(\frac{\imath}{2})^{n-4}(-1)^{\frac{(n-2)(n-3)}{2}}(X_1\wedge\dots\wedge X_{n-3}\wedge X)\wedge
(\overline{X}_1\wedge\dots\wedge\overline{X}_{n-3} \wedge\overline{X}).\\
\end{array}$$

\noindent Simple computations give

\begin{eqnarray*}
X_1\wedge\dots\wedge X_{n-3}\wedge X & = & z^{n-4}_{n-2}(|z_1|^2+\dots+|z_{n-3}|^2-|z_{n-2}|^2)dz_1\dots dz_{n-2}\\
& = & - z^{n-4}_{n-2}dz_1\dots dz_{n-2}.\\
\end{eqnarray*}

\noindent Similarly, $\displaystyle{ \overline{X}_1\wedge\dots\wedge\overline{X}_{n-3}\wedge\overline{X}=-\bar{z}^{n-4}_{n-2}d\bar{z}_1\dots
d\bar{z}_{n-2}}$. Therefore,

$$\begin{array}{lcl}

(X_1\wedge\dots\wedge X_{n-3}\wedge X)\wedge (\overline{X}_1\wedge\dots\wedge\overline{X}_{n-3}\wedge\overline{X})&=&
|z_{n-2}|^{2(n-4)}dz_1\dots dz_{n-2} d\bar{z}_1\dots d\bar{z}_{n-2}\\
&=&2^{n-2}\imath^{(n-2)(n-4)}|z_{n-2}|^{2(n-4)}d\lambda_{2(n-2)},\\
\end{array}$$

\noindent and we get

$$U^*_1\wedge V^*_1\wedge\dots\wedge U^*_{n-3}\wedge V^*_{n-3}\wedge df \wedge (df\circ\J)=4|z_{n-2}|^{2(n-4)}d\lambda_{2(n-2)}.$$

\noindent By the definition of $\hat{\mu}^1_\mathrm{Tr}$, we obtain

$$\hat{\mu}^1_\mathrm{Tr}(\S^1\cdot Z)=\frac{c_0}{4|z_{n-2}|^{2(n-4)}}U^*_1\wedge V^*_1\wedge\dots\wedge U^*_{n-3}\wedge V^*_{n-3}.$$

\carre

\rem \begin{itemize}

\item[-] Even though the $1$-forms $U^*_k$ and $V^*_k$ are not invariant by the $\S^1$-action, the $2$-form $U^*_k\wedge V^*_k$ is. Hence, the
$2(n-3)$-form $U^*_1\wedge V^*_1\wedge \dots \wedge U^*_{n-3}\wedge V^*_{n-3}$ is invariant by the $\S^1$ action, and we get a well defined
volume form on $\mathbf{Q}_1/\S^1$.

\item[-] Let $\mu^1_\mathrm{Tr}$ be the volume form on $\mathbf{Q}_1$ verifying the following condition

$$ \mu^1_\mathrm{Tr}\wedge df =\mu_\mathrm{Tr}.$$

The tangent vector to the $\S^1$ orbit at a point $Z \in \C^2$ is given by $\imath Z$, and we have

$$df\circ\J (\imath Z)=-df(Z)=-\langle Z,Z \rangle =1.$$

Therefore, the volume form $\hat{\mu}^1_\mathrm{Tr}$ can be considered as the push-forward of $\mu^1_\mathrm{Tr}$ onto $\mathbf{Q}_1/\S^1$.

\end{itemize}

Now, we will proceed to compute the volume form defined by $\eta$ on $Z^\bot$ in terms of $U^*_k,V^*_k$. Let $(\eta_{ij})$ with
$i,j=1,\dots,2(n-3)$ be the (real) matrix of $\eta$ in the basis $\{U_1,V_1,\dots,U_{n-3},V_{n-3}\}$. Since the volume form $\mu_\mathrm{Hyp}$
is defined by the metric $\eta$, we have

\begin{equation}\label{VolHyp}
\mu_\mathrm{Hyp}(\S^1\cdot Z)=\frac{1}{\sqrt{\det(\eta_{ij})}}U^*_1\wedge V^*_1\wedge\dots\wedge U^*_{n-3}\wedge V^*_{n-3}
\end{equation}

\begin{lemma}\label{SSlmC2}

We have $\displaystyle{\det(\eta_{ij})=|z_{n-2}|^{4(n-4)}}$.

\end{lemma}

\dem Since $\eta$ is the real part of $\Her$, the matrix $(\eta_{ij})$ is the real interpretation of the matrix $(\Her_{ij})$ of $\Her$ in the
complex basis $\{U_1,\dots,U_{n-3}\}$ of $Z^\bot$. This implies $\displaystyle{\det(\eta_{ij})=|\det(\Her_{ij})|^2}$.\\

\noindent We have $\displaystyle{\Her_{ij}=\langle U_i,U_j\rangle=\left\{%
\begin{array}{ll}
    -z_i\bar{z}_j, & \hbox{  if  $i \neq j$;} \\
    |z_{n-2}|^2-|z_i|^2, & \hbox{ if  $i=j$.} \\
\end{array}%
\right.}$, hence

$$\det(\Her_{ij})= |z_{n-2}|^{2(n-3)} \left|%
\begin{array}{cccc}
  1-|\varepsilon_1|^2 & -\bar{\varepsilon}_1\varepsilon_2 & \dots & -\bar{\varepsilon}_1\varepsilon_{n-3} \\
  -\bar{\varepsilon}_2\varepsilon_1 & 1-|\varepsilon_2|^2 & \dots & -\bar{\varepsilon}_2\varepsilon_{n-3} \\
  \cdots & \cdots & \cdots & \cdots \\
  -\bar{\varepsilon}_{n-3}\varepsilon_1 & -\bar{\varepsilon}_{n-3}\varepsilon_2 & \dots & 1-|\varepsilon_{n-3}|^2 \\
\end{array}%
\right|, $$

\noindent where $\varepsilon_k=z_k/z_{n-2},\; k=1,\dots,n-3$. Since we have

$$\left|%
\begin{array}{cccc}
  1-|\varepsilon_1|^2 & -\bar{\varepsilon}_1\varepsilon_2 & \dots & -\bar{\varepsilon}_1\varepsilon_{n-3} \\
  -\bar{\varepsilon}_2\varepsilon_1 & 1-|\varepsilon_2|^2 & \dots & -\bar{\varepsilon}_2\varepsilon_{n-3} \\
  \cdots & \cdots & \cdots & \cdots \\
  -\bar{\varepsilon}_{n-3}\varepsilon_1 & -\bar{\varepsilon}_{n-3}\varepsilon_2 & \dots & 1-|\varepsilon_{n-3}|^2 \\
\end{array}%
\right|= 1-(|\varepsilon_1|^2+\dots+|\varepsilon_{n-3}|^2),$$

\noindent it follows that

\begin{eqnarray*}
\det(\Her_{ij})&=&|z_{n-2}|^{2(n-3)}(1-(|\varepsilon_1|^2+\dots+|\varepsilon_{n-3}|^2))\\
&=&|z_{n-2}|^{2(n-4)}(|z_{n-2}|^2-(|z_1|^2+\dots+|z_{n-3}|^2))\\
&=& |z_{n-2}|^{2(n-4)}.\\
\end{eqnarray*}

\noindent Consequently, we have $\det(\eta_{ij})=|\det(\Her_{ij})|^2=|z_{n-2}|^{4(n-4)}$. The lemma is then proved. \carre

\subsection{Proof of Proposition \ref{SSprC}}

From Lemma \ref{SSlmC1}, Lemma \ref{SSlmC2}, and (\ref{VolHyp}), we know that the quotient $\displaystyle{
\frac{\hat{\mu}^1_\mathrm{Tr}}{\mu_\mathrm{Hyp}} }$ is a locally constant function on $\MSSvi$. It remains to show that $\MSSvi$ is connected,
but this is known since, by Proposition 2 \cite{Troy}, we can identify $\MSSvi$ to modular space of $n$-punctured sphere, which is
connected.\carre

\begin{appendices}

\section{Proof of Theorem \ref{ETthF}}\label{prfFSect}

Let $n_1$ be the cardinal of  $V\cap\partial\Sig$, and $n_2$ be the cardinal of $V\cap\inter(\Sig)$. Observe that any triangulations of $\Sig$
whose vertex set is $V$ has a fixed number $N_e$ of edges. Let $k, 0\leq k \leq N_e$, be the number of common edges of $\T_1$ and $\T_2$. Since
the boundary of $\Sig$ contains $n_1$ edges, we have $k\geq n_1$. Assume that $n_1 \leq k <N_e$, we will proceed by induction.\\

First, let us prove the following technical lemma

\begin{lemma}\label{ETlmF2}
Let $\P$ be a polygon in $\R^2$ whose vertices are denoted by $A_1,A_2,A_3, B_1,\dots,B_l$. Let $x:\R^2\lra \R$, and $y: \R^2\lra \R$ denote the
two coordinate functions of $\R^2$. Assume that the vertices of $\P$ verify the following conditions:

\begin{itemize}
\item[.] $(A_1,A_2,A_3)$ are ordered in the clockwise sense.

\item[.] $y(A_i)\geq 0,\; i=1,2,3$,  $y(A_1)<y(A_2)$, and $y(A_2)\geq y(A_3)$.

\item[.] $y(B_j)<0, \; j=1,\dots,l$;

\item[.] $B_1,\dots,B_l$ are ordered in the counter-clockwise sense.

\item[.] For all $j\in\{1,\dots,l\}$, the segment $\overline{A_2B_j}$ is a diagonal of $\P$.

\end{itemize}

\noindent Let $\T$ denote the triangulation of $\P$ by the diagonals $\overline{A_2B_1},\dots,\overline{A_2B_l}$. Let $\{s_0,\dots,s_k\}$ be a
family of disjoint horizontal segments in $\P$ whose endpoints are contained the boundary of $\P$, where $s_0$ is a segment lying on the
horizontal axis $y=0$. Let $r$ be the number of intersection points of the edges of $\T$ with the set $\DS{\cup_{i=0}^k s_i}$. Then there exists
a sequence of elementary moves which transform $\T$ into a new triangulation $\T'$ whose edges intersect the set $\DS{\cup_{i=0}^k s_i}$ at at
most $r-1$ points.

\end{lemma}

\dem Let $j_0$ be the smallest index such that $y(B_{j_0})=\min\{y(B_j) : j=1,\dots,l\}$, and consider the following algorithm:

\begin{enumerate}

\item If $\P$ is a quadrilateral, that is $l=1$, then $\P$ must be convex since its two diagonals intersect. Apply an elementary move inside
$\P$ and stop the algorithm.

\item If $1<j_0<l$, then consider the quadrilateral $A_{2}B_{j_0-1}B_{j_0}B_{j_0+1}$. By the choice of $j_0$, this quadrilateral is convex,
hence, we can apply an elementary move inside it, and  the algorithm stops.

\item If $j_0=1$ and $l\geq 2$, then consider the quadrilateral $A_{2}A_{1}B_{1}B_{2}$. Observe that this quadrilateral is convex. Apply an
elementary move inside it. By this move, we get a new triangulation of $\P$ which contains the triangle $\Delta A_{1}B_{1}B_{2}$. Cut off this
triangle from $\P$. Replace $\P$ by the remaining sub-polygon and restart the algorithm.

\item  If $j_0=l>1$, then consider the quadrilateral $A_{2}A_{3}B_{l}B_{l-1}$. Since this quadrilateral is convex, we can apply an elementary
move inside it, then cut off the triangle $\Delta A_{3}B_{l}B_{l-1}$. Replace $\P$ by the remaining sub-polygon and restart the algorithm.

\end{enumerate}

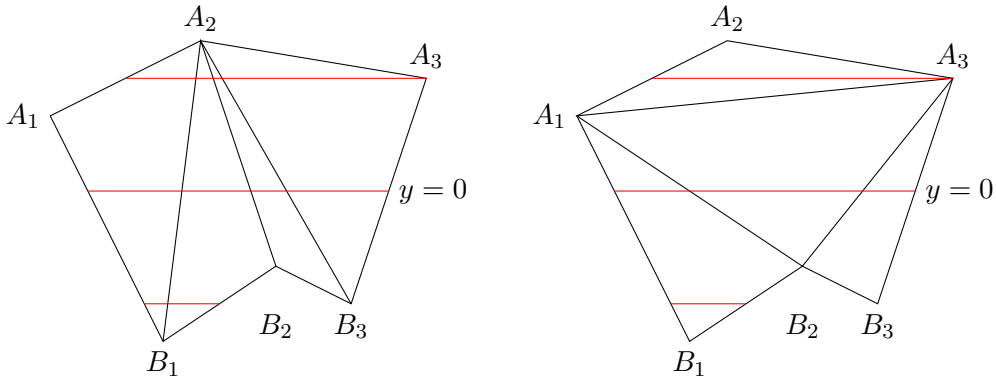
\begin{center}
\begin{tikzpicture}
\draw (-0.5,1) -- (1.5,2) -- (4.5,1.5) -- (3.5,-1.5) -- (2.5,-1) -- (1,-2) -- cycle;  \draw (6.5,1) -- (8.5,2) -- (11.5,1.5) -- (10.5,-1.5) --
(9.5,-1) -- (8,-2) -- cycle;

\draw (1,-2) -- (1.5,2) -- (2.5,-1); \draw (1.5,2) -- (3.5,-1.5); \draw (6.5,1) -- (9.5,-1) -- (11.5,1.5) -- cycle;

\draw (-0.5,1) node[left] {$A_{1}$}; \draw (1.5,2) node[above] {$A_{2}$}; \draw (4.5,1.5) node[above] {$A_{3}$}; \draw (1,-2) node[below]
{$B_{1}$}; \draw (2.5,-1.5) node[below] {$B_{2}$}; \draw (3.5,-1.5) node[below] {$B_{3}$};

\draw (6.5,1) node[left] {$A_{1}$}; \draw (8.5,2) node[above] {$A_{2}$}; \draw (11.5,1.5) node[above] {$A_{3}$}; \draw (8,-2) node[below]
{$B_{1}$}; \draw (9.5,-1.5) node[below] {$B_{2}$}; \draw (10.5,-1.5) node[below] {$B_{3}$};

\draw[red] (0,0) -- (4,0); \draw[red] (0.5,1.5) -- (4.5,1.5); \draw[red] (0.75,-1.5) -- (1.75,-1.5); \draw[red] (7,0) -- (11,0); \draw[red]
(7.5,1.5) -- (11.5,1.5); \draw[red] (7.75,-1.5) -- (8.75,-1.5);

\draw (4,0) node[right] {$y=0$}; \draw (11,0) node[right] {$y=0$};

\end{tikzpicture}
\end{center}

\noindent Observe that, at each step of the algorithm above, the number of intersection points of the set $\cup_{i=0}^k s_i$ with the edges of
the new triangulation cannot exceed the number of intersection points with those of the ancien one. Indeed, suppose that we are in the case
$2.$, by the choice of $j_0$, we have $y(B_{j_0})\leq \min\{y(B_{j_0-1}),y(B_{j_0+1})\}$, and $y(A_{2})\geq \max\{y(B_{j_0-1}),y(B_{j_0+1})\}$,
consequently, if a horizontal segment $s_i$ intersects $\overline{B_{j_0-1}B_{j_0+1}}$, then it must intersect $\overline{A_{2}B_{j_0}}$.
Therefore, the number of intersection points does not increase. The same argument works for the other cases.\\

\noindent Moreover, at the final step of the algorithm, \ie case 1. or 2.,  we replace a diagonal intersecting $s_0$ by another one which does
not intersect $s_0$. Therefore, by this algorithm, we get a new triangulation $\T'$ of $\P$ whose edges have strictly less intersection points
with the set $\cup_{i=0}^k s_i$ than those of $\T_4$.\carre

\vspace{0.5cm}

Let $a_1,\dots,a_{N_e}$, and $b_1,\dots,b_{N_e}$ denote the edges of $\T_1$ and $\T_2$ respectively. We can assume that $a_i=b_i$, for
$i=1,\dots,k$.  All we need to prove is the following

\begin{proposition}\label{ETprF1}

There exists a sequence of elementary moves which transform $\T_1$ into a new triangulation containing $b_1,\dots,b_k$, and $b_{k+1}$.

\end{proposition}

\dem Since $b_{k+1}$ is not an edge of $\T_1$, it must intersect some edges of $\T_1$. Let $\P$ be the developing polygon of $b_{k+1}$ with
respect to $\T_1$, and $\varphi: \P \lra \Sig$ be the associated immersion. Let $\T_3$ be the triangulation of $\P$ by diagonals which is
induced by $\T_1$, (\ie $\T_3=\varphi^{-1}(\T_1)$). By definition, each diagonal in $\T_3$ is mapped by $\varphi$ onto an edge of $\T_1$ which
intersects $\inter(b_{k+1})$. Let $d$ be the diagonal of $\P$ such that $\varphi(d)=b_{k+1}$. Observe that $\inter(d)$ intersects all the
diagonals which are edges of $\T_3$.\\

\noindent Let $m$ be the number of intersection points of $b_{k+1}$ with the edges of $\T_1$ excluding the two endpoints of $b_{k+1}$. Note that
$b_{k+1}$ may intersect an edge of $\T_1$ more than once. By construction, the polygon $\P$ is triangulated by $m$ diagonals, hence it has $m+3$
sides. We prove the proposition by induction.

\begin{itemize}
\item[-] If $m=1$, then $\P$ is a quadrilateral. The quadrilateral $\P$ must be convex because its two diagonals intersect. If $\P$ is mapped by
$\varphi$ to a single triangle of $\T_1$, then there is a singular point of $\Sig$ with cone angle strictly less than $\pi$, but this is
impossible since we have assumed that $\Sig$ verifies Property $\pQ$. Thus, we conclude that $\varphi$ maps $\inter(\P)$ isometrically onto a
domain consisting of two triangles in $\T_1$. Clearly, by applying the elementary move inside $\varphi(\P)$, we obtains a new triangulation
which contains $b_{k+1}$.\\

\item[-] If $m>1$, it is enough to show that there exists a sequence of elementary moves which transform $\T_1$ into a new triangulation $\T'_1$
containing $b_1=(a_1),\dots,b_k=(a_k)$, such that $b_{k+1}$ intersects the edges of $\T'_1$ at most $m-1$ times.\\

Equip the plane $\R^2$ with a system of Cartesian coordinates such that the diagonal $d \subset \P$ is a horizontal segment lying in the axis
$Ox$. Let $x:\R^2\lra \R$, and $y:\R^2\lra \R$ denote the two coordinate functions. Let $A_1,\dots,A_r$ denote the vertices of $\P$ such that
$y(A_i)>0$, and $B_1,\dots,B_s$ denote the vertices of $\P$ such that $y(B_j)<0$. Let $A_0$ and $A_{r+1}$ denote the left and the right
endpoints of $d$ respectively. We set, by convention, $B_0=A_0$, and $B_{s+1}=A_{r+1}$. Since $\P$ has $m+3$ vertices, we have $r+s+2=m+3$. We
can assume that $r \geq s$ (if it is not the case, reverse the orientation of $Oy$). We denote the vertices of $\P$ so that $A_0,\dots,A_{r+1}$
are ordered in the clockwise sense, and $B_0,\dots,B_{s+1}$ are ordered in the counter-clockwise sense.\\

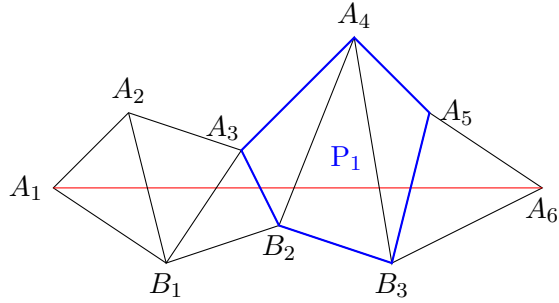
\begin{center}
\begin{tikzpicture}
\draw (0,0) -- (1,1) -- (2.5,0.5) -- (4,2) -- (5,1) -- (6.5,0) -- (4.5,-1) -- (3,-0.5) -- (1.5,-1)  -- cycle;

\draw[red] (0,0) -- (6.5,0); \draw (1,1) -- (1.5,-1) -- (2.5,0.5) -- (3,-0.5) -- (4,2) -- (4.5,-1) -- (5,1); 
\draw (0,0) node[left] {$A_1$}; \draw (1,1) node[above] {$A_2$}; \draw (2.25,0.5) node[above] {$A_3$}; \draw (4,2) node[above] {$A_4$}; \draw
(5,1) node[right] {$A_5$}; \draw (6.5,0) node[below] {$A_6$}; \draw (1.5,-1) node[below] {$B_1$}; \draw (3,-0.5) node[below] {$B_2$}; \draw
(4.5,-1) node[below] {$B_3$};

\draw[thick, blue] (2.5,0.5) -- (4,2) -- (5,1) -- (4.5,-1) -- (3,-0.5) -- cycle;

\draw[blue] (3.9,0.1) node[above] {$\P_1$};

\end{tikzpicture}
\end{center}

Without loss of generality, we can assume that $r\geq 2$ because $m>1$. Let $i_0$ be the smallest index such that $y(A_{i_0})=\max\{y(A_i):
i=1,\dots,r\}$, that is $y(A_{i_0})\geq y(A_i) \; \forall i=1,\dots,r$, and $y(A_{i_0})>y(A_i)$ if $i<i_0$. Consider the polygon $\P_1$ which is
the union of all triangles in $\T_3$ having $A_{i_0}$ as a vertex. There exist $j_0$ and $l$ such that the vertices of $\P_1$ are $A_{i_0-1},
A_{i_0}, A_{i_0+1}$ and $B_{j_0},\dots,B_{j_0+l}$, note that  $\P_1$ is triangulated by the diagonals
$\overline{A_{i_0}B_{j_0}},\dots,\overline{A_{i_0}B_{j_0+l}}$. Let $\T_4$ denote this triangulation of $\P_1$.\\

By Lemma \ref{ETlmF3} below, we know that the restriction of $\varphi$ into  $\inter(\P_1)$ is injective. Let $\mathrm{Q}_1$ be the image of
$\inter(\P_1)$ under $\varphi$. Since $b_1,\dots,b_k,b_{k+1}$ are edges of the triangulation $\T_2$, we have
$\inter(b_i)\cap\inter(b_{k+1})=\vide, \; \forall i=1,\dots,k$. Recall that $b_1,\dots,b_k$ are also edges of the triangulation $\T_1$ of
$\Sig$, from this we deduce that $\inter(b_i)\cap\mathrm{Q}_1=\vide$, since if $e$ is an edge of $\T_1$ and $\inter(e)\cap \mathrm{Q}_1\neq
\vide$, then $\inter(e)\cap\inter(b_{k+1})\neq\vide$. This implies that an elementary move inside $\mathrm{Q}_1$ does not affect the edges
$b_1,\dots,b_k$.\\

Consider the intersection of  $\P_1$ and  $\varphi^{-1}(b_{k+1})$. A priori, this set is a family of geodesic segments with endpoints in the
boundary of $\P_1$. Clearly, the segment $s_0=\overline{A_0A_{r+1}}\cap\P_1$ is contained in the set $\P_1\cap \varphi^{-1}(b_{k+1})$. Since
$\Sig$ satisfies $\pQ$, all the segments in this family are parallel, therefore, all of them are parallel to the horizontal axis. Let $\delta$
be the number of intersection points of the set $\P_1\cap\varphi^{-1}(b_{k+1})$ and the edges of $\T_4$.\\

Now, Lemma \ref{ETlmF2} shows that there exists a sequence of elementary moves which transform $\T_4$ into a new triangulation of $\P_1$ whose
edges intersect the set $\P_1\cap\varphi^{-1}(b_{k+1})$ at at most $\delta-1$ points. It follows that there exists a sequence of elementary
moves inside the domain $\mathrm{Q}_1 \subset \Sig$ which transform $\T_1$ into a new triangulation of $\Sig$ whose edges have at most $m-1$
intersection points with $b_{k+1}$. As we have seen, those elementary moves do not affect the edges $b_1,\dots,b_k$. By induction, the
proposition is then proved. \carre

\end{itemize}

\noindent We need the following lemma to complete the proof of Proposition \ref{ETprF1}

\begin{lemma}\label{ETlmF3}

With the same notations as in the proof of \ref{ETprF1}, the restriction of $\varphi$ to $\inter(\P_1)$ is an isometric embedding.

\end{lemma}

\dem Since $\varphi$ maps each triangle of $\T_3$ onto a triangle of $\T_1$, it is enough to show that the images by $\varphi$ of the triangles
of $\T_3$ which are contained in $\P_1$ are all distinct. Suppose that there exist two triangles $\Delta_1$ and $\Delta_2$ such that
$\varphi(\Delta_1)=\varphi(\Delta_2)$. Since $\varphi$ is locally isometric, and by assumption, the orthogonal part of the holonomy of any
closed curve in $\inter(\Sig\setm\{p_1,\dots,p_n\})$ is either $\Id$ or $-\Id$, it follows that either $\Delta_2=\Delta_1+v$, or
$\Delta_2=-\Delta_1+v$, where $-\Delta_1$ is the image of $\Delta_1$ by $-\Id$, and $v\in \R^2$. Note that, by definition, the triangles
$\Delta_1$ and $\Delta_2$ have a common vertex, which is $A_{i_0}$.

\begin{itemize}

\item[$\bullet$] If $\Delta_2=\Delta_1+v$, excluding the case $\Delta_1\equiv\Delta_2$, we have two possible configurations. In these both
cases, we see that the angle of $\P_1$ at the point $A_{i_0}$ is at least $\pi$. But, by assumption, this is impossible since we have
$y(A_{i_0})>y(A_{i_0-1})$ and $y(A_{i_0})\geq y(A_{i_0+1})$.

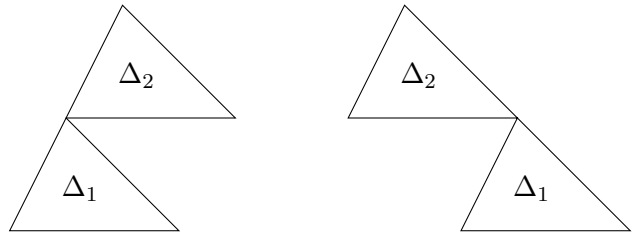
\begin{center}
\begin{tikzpicture}[scale=0.75]
\draw (0,0) -- (2,-2) -- (-1,-2) -- cycle; \draw (0,0) -- (3,0) -- (1,2) -- cycle;   \draw (0.25,-1.25) node {$\Delta_1$}; \draw (1.25,0.75)
node {$\Delta_2$};

\draw (5,0) -- (6,2) -- (10,-2) -- (7,-2) -- (8,0) -- cycle; \draw (6.25,0.75) node {$\Delta_2$}; \draw (8.25,-1.25) node {$\Delta_1$};

\end{tikzpicture}
\end{center}

\item[$\bullet$] If $\Delta_2=-\Delta_1+v$, we have three possible configurations. In the case where $\Delta_1$ and $\Delta_2$ have only one
common vertex, we see that the angle of $\P_1$ at $A_{i_0}$ must be greater than $\pi$, which is, as we have seen above, impossible. In the
other two cases, $\Delta_1$ and $\Delta_2$ are adjacent. As we have seen,  this implies the existence of a singular point of $\Sig$ with cone
angle strictly less than $\pi$. This is again impossible.

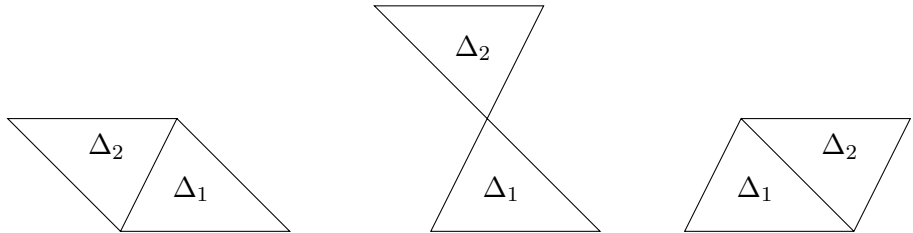
\begin{center}
\begin{tikzpicture}[scale=0.75]
\draw (0,0) -- (2,-2) -- (-1,-2) -- cycle; \draw (0.25,-1.25) node {$\Delta_1$}; \draw (0,0) -- (1,2) -- (-2,2) -- cycle;  \draw (-0.25,1.25)
node {$\Delta_2$};

\draw (-5.5,0) -- (-3.5,-2) -- (-6.5,-2) -- (-5.5,0) -- (-8.5,0) -- (-6.5,-2); \draw (-5.25,-1.25) node {$\Delta_1$}; \draw (-6.75,-0.5) node
{$\Delta_2$};

\draw (4.5,0) -- (7.5,0) -- (6.5,-2) -- (4.5,0) -- (3.5,-2) -- (6.5,-2); \draw (4.75,-1.25) node {$\Delta_1$}; \draw (6.25,-0.5) node
{$\Delta_2$};


\end{tikzpicture}
\end{center}

\end{itemize}

\noindent The lemma is then proved. \carre

\section{ Proof of Theorem \ref{SSThD}}\label{prfSSThD}

Let $x_1,\dots,x_n$ denote the vertices of $\T_1$ and $\T_2$. By convention, we consider $\{x_1,\dots,x_n\}$ as the set of singular points of
$\Sig$ even though some of them may be regular. In what follows, if $\T$ is a triangulation of $\Sig$ whose vertex set is $\{x_1,\dots,x_n\}$,
we will call a tree contained in the $1$-skeleton of $\T$ which connects all the vertices of $\T$ a {\em maximal tree}. Let $A_i, \; i=1,2$ be a
maximal tree of $\T_i$. If $A_1\equiv A_2$, then the theorem follows from Theorem \ref{ETthF}. Thus, it is enough to prove the following

\begin{proposition}\label{SSprD1}
There exists a sequence of elementary moves which transforms $\T_1$ into a triangulation containing $A_2$.
\end{proposition}

\noindent We start by the following lemma

\begin{lemma}\label{SSlmD1}
If $c_1,\dots,c_k$ are geodesic segments with endpoints in $\{x_1,\dots,x_n\}$ such that $\inter(c_i)\cap\inter(c_j)=\vide$ if $i\neq j$, and
$\inter(c_i)\cap A_1=\vide, i=1,\dots,k$, then there exists a sequence of elementary moves which transforms $\T_1$ into a new
triangulation containing $A_1$, and all the segments $c_1,\dots,c_k$. \\

\end{lemma}

\dem Let $\P_{c_1}$ be the developing polygon of $c_1$ with respect to $\T$. Since $\inter(c_1)\cap A_1=\vide$, the isometric immersion
$\displaystyle{\varphi_{c_1}: \inter(\P_{c_1})\ra \Sig}$ is an embedding, therefore, by applying Theorem \ref{ETthF} to $\P_{c_1}$, we deduce
that $\T$ can be transformed by elementary moves into a triangulation containing $A_1$ and $c_1$. We can then restart the procedure with $c_2$
in the place of $c_1$. Since $\inter(c_1)\cap\inter(c_2)=\vide$, the developing polygon of $c_2$ does not contain $c_1$ in the interior. Thus,
the lemma follows by induction.\carre

Now, let $a_1,\dots,a_{n-1}$ denote the edges of the tree $A_1$, and $b_1,\dots,b_{n-1}$ denote the edges of the tree $A_2$. We will proceed by
induction. Suppose that $\T_1$ contains already the $k$ edges $b_1,\dots,b_k$ of $A_2$.  We will show that $\T_1$ can be transformed by a
sequence of elementary moves into a new triangulation containing $b_1,\dots,b_k$ and $b_{k+1}$. Let $m$ be the number of intersection points of
$b_{k+1}$ with the tree $A_1$ excluding the endpoints of $b_{k+1}$. If $m=0$, then Lemma \ref{SSlmD1} allows us to get the conclusion.
Therefore, if $m\geq 1$, all we need to show is the following

\begin{lemma}\label{SSlmD2}
The triangulation $\T_1$ can be transformed by elementary moves into a new triangulation $\T'_1$ which contains a maximal tree $A'_1$, and the
edges $b_1,\dots,b_k$, such that the number of intersecting points of $b_{k+1}$ with $A'_1$, excluding the endpoints of $b_{k+1}$, is at most
$m-1$.
\end{lemma}

\dem We can assume that the endpoints of $b_{k+1}$ are $x_1$ and $x_2$. We consider $b_{k+1}$ as a geodesic ray exiting from $x_1$. Let $y_1$
denote the first intersection point of $b_{k+1}$ with the tree $A_1$, which  is contained in the interior of an edge
$\overline{x_{j_1}x_{j_1+1}}$ of  $A_1$. Let $\overline{x_1y_1}$ denote the subsegment of $b_{k+1}$ whose endpoints are $x_1$ and $y_1$. Without
loss of generality, we can assume that $x_{j_1}$ is contained in the unique path along $A_1$ from $x_1$ to $x_{j_1+1}$.\\

\noindent Cutting open the surface $\Sig$ along the tree $A_1$, we  get a flat surface $\Sig'$ with geodesic boundary homeomorphic to a close
disk. By construction, we have a surjective map $\DS{\pi_{A_1} : \Sig'\lra \Sig}$ verifying the following properties

\begin{itemize}
\item[.] $\pi_{A_1}|_{\inter(\Sig')}$ is an isometry,

\item[.] $\pi_{A_1}(\partial \Sig')=A_1$.

\item[.] There are $2(n-1)$ geodesic segments in the boundary  of $\Sig'$ such that the restriction of $\pi_{A_1}$ into each of which is an
isometry.

\item[.] For every edge $e$ in $A_1$, $\pi_{A_1}^{-1}(\inter(e))$ is the union of two open segments in the boundary of $\Sig'$.

\end{itemize}

\noindent Let $s_1$ denote $\pi_{A_1}^{-1}(\overline{x_1y_1})$, then $s_1$ is a geodesic segment with endpoints in $\partial\Sig'$. Let $x'_1$
and $y'_1$ denote the endpoints of $s_1$ with $\pi_{A_1}(x'_1)=x_1$, and $\pi_{A_1}(y'_1)=y_1$. Let $x'_1,\dots,x'_{2(n-1)}$ denote the points
in $\pi^{-1}_{A_1}(\{x_1,\dots,x_n\})$ following an orientation of $\partial \Sig'$. By choosing the suitable orientation, we can assume that
the point $y'_1$ is between $x'_{j'_1}$ and $x'_{j'_1+1}$, where $\pi_{A_1}(x'_{j'_1})=x_{j_1}$, and $\pi_{A_1}(x'_{j'_1+1})=x_{j_1+1}$. For
every $j$ in $\{1,\dots,2(n-1)\}$, we denote by $\overline{x'_jx'_{j+1}}$ the segment in the boundary of $\Sig'$ with endpoints $x'_j$ and
$x'_{j+1}$, with the convention $x'_{2n-1}=x'_1$. Note that $\pi_{A_1}(\overline{x'_jx'_{j+1}})$ is an edge of $A_1$.\\

\noindent Let $c_0$ be a path in $\Sig'$ joining $x'_1$ and $x'_{j'_1+1}$ with minimal length. First, we prove

\begin{lemma}\label{SSlmD3}
 We have $c_0\cap s_1=\{x'_1\}$.

\end{lemma}

\dem Suppose that $c_0\cap\inter(s_1)\neq\vide$, then let $y'_2$ denote the first intersection point of $c_0$ with $s_1\setm\{x'_1\}$. Let $c_1$
denote the path from $x'_1$ to $y'_2$ along $c_0$, and let $\overline{x'_1y'_2}$ denote the  subsegment of $s_1$ with endpoints $x'_1$ and
$y'_2$. By definition, we see that $c_1\cup\overline{x'_1y'_2}$ is the boundary a flat disk with piecewise geodesic boundary $\mathrm{D}$. Since
the path $c_0$ is of minimal length,  it follows that the interior angle between two consecutive segments of $c_1$ is at least $\pi$. Therefore,
if the number of segments in $c_1$ is $l$, the boundary of $\mathrm{D}$ contains then  $l+1$ geodesic segments, and the sum of all the interior
angles is at least $(l-1)\pi$. But this is impossible by the Gauss-Bonnet Theorem, hence we conclude that $c_0\cap(s_1\setm\{x'_1\})=\vide$.
\carre

\noindent Let $\overline{y'_1x'_{j'_1+1}}$ denote the subsegment of $\overline{x'_{j'_1}x'_{j'_1+1}}$ with endpoints $y'_1$ and $x'_{j'_1+1}$.
From Lemma \ref{SSlmD3}, we see that $s_1\cup\overline{y'_1x'_{j'_1+1}}\cup c_0$ is the boundary of a disk $\mathrm{D}_0$ contained in $\Sig'$.
We have immediately the following

\begin{lemma}\label{SSlmD4}
Let $s$ be a geodesic ray that intersects the interior of $\mathrm{D}_0$. If $s$ inters $\mathrm{D}_0$ by a point in the path $c_0$, then $s$
must exit $\mathrm{D}_0$ by  a point in $(s_1\cup\overline{y'_1x'_{j'_1+1}})\setm\{x'_1,x'_{j'_1+1}\}$.

\end{lemma}

\dem If $s$ exits $\mathrm{D}_0$ by another point in $c_0$, then we have a flat disk with piecewise geodesic boundary which violates the
Gauss-Bonnet Theorem.\carre

\noindent Let $\hat{c}_0$ denote the image of $c_0$ under $\pi_{A_1}$. The path $\hat{c}_0$ is then a finite union of geodesic segments on
$\Sig$ with endpoints in the set $\{x_1,\dots,x_n\}$. It is clear that $\hat{c}_0$ contains a path $\hat{c}_1$ joining $x_1$ to $x_{j_1+1}$. Let
us prove the following

\begin{lemma}\label{SSlmD5}
 The path $\hat{c}_1$ does not contain the segment $\overline{x_{j_1}x_{j_1+1}}$.
\end{lemma}

\dem Suppose, on the contrary, that $\hat{c}_1$ contains $\overline{x_{j_1}x_{j_1+1}}$. This implies that $c_0$ contains a segment
$\overline{x'_{k'}x'_{k'+1}}$, with $k'\neq j'_1$, such that
$\DS{\pi_{A_1}(\overline{x'_{k'}x'_{k'+1}})=\pi_{A_1}(\overline{x'_{j'_1}x'_{j'_1+1}})= \overline{x_{j_1}x_{j_1+1}}}$.\\

\noindent Let $y'_2$ denote the unique point in $\overline{x'_{k'}x'_{k'+1}}$ such that $\pi_{A_1}(y'_2)=\pi_{A_1}(y'_1)=y_1$. Observe that
$\pi^{-1}_{A_1}(b_{k+1})$ is a sequence of $(m+1)$ geodesic segments of $\Sig'$ with endpoints in the boundary of $\Sig'$, $s_1$ is the first
one of this sequence, $y'_2$ is then one endpoint of  the next segment in this sequence, which will be denoted by $s_2$.\\

\noindent  By assumption, $y'_2$ is an intersection point of the segment $s_2$ and the disk $\mathrm{D}_0$. Consider the segment $s_2$ as a
geodesic ray exiting from $y'_2$. By Lemma \ref{SSlmD4}, the ray $s_2$ exits $\mathrm{D}_0$ by a point $z'_2$ in
$(s_1\cup\overline{y'_1x'_{j'_1+1}})\setm\{x'_1,x'_{j'_1+1}\}$. Since  the geodesic $b_{k+1}$ is a simple, $z'_2$ can not be a point in $s_1$,
hence $z'_2$ must be a point in $\inter(\overline{y'_1x'_{j'_1+1}})$.\\

\noindent Now, since the segments $\overline{x'_{j'_1}x'_{j'_1+1}}$ and $\overline{x'_{k'}x'_{k'+1}}$ are identified by $\pi_{A_1}$, the point
$z'_2$ is identified to a point $y'_3$ in $\overline{x'_{k'}x'_{k'+1}}$. Consequently, the argument above can be applied infinitely many times,
which implies that $\pi^{-1}_{A_1}(b_{k+1})$ contains infinitely many segments, and we have a contradiction to the fact that
$\pi^{-1}_{A_1}(b_{k+1})$ contains only $m+1$ segments.  \carre

Since $A_1$ is a tree, the set $A_1\setm\inter(\overline{x_{j_1}x_{j_1+1}})$ has two connected components, the one containing $x_1$ will be
denoted by $C_1$, the other one containing $x_{j_1+1}$ will be denoted by $C_2$. From Lemma \ref{SSlmD5}, we know that the path $\hat{c}_1$,
which joins $x_1$ to $x_{j_1+1}$ does not contain $\overline{x_{j_1}x_{j_1+1}}$. Therefore the path $\hat{c}_1$ must contain a segment
$\hat{s}$, with endpoints in $\{x_1,\dots,x_n\}$, such that one of the two endpoints is in $C_1$, and the other is in $C_2$. Let $s$ denote
$\pi^{-1}_{A_1}(\hat{s})$. Evidently, $\hat{s}$ is not an edge of $A_1$, hence $s$ is a segment contained inside $\Sig'$, it follows that
$\inter(\hat{s})\cap A_1=\vide$. Let us prove

\begin{lemma}\label{SSlmD6}
$\inter(\hat{s})\cap\inter(b_i)=\vide$, for every $i=1,\dots,k$.
\end{lemma}

\dem Let $b'_i, \; i=1,\dots,k$, denote $\pi^{-1}_{A_1}(b_i)$. Since $\inter(b_i)\cap A_1=\vide$, $b'_i$ is a geodesic segment with endpoints in
$\partial \Sig'$.  Suppose that $\inter(\hat{s})\cap\inter(b_i) \neq \vide$, it follows that $\inter(b'_i)\cap\inter(s)\neq \vide$. Let $y''_i$
be the intersection point of $\inter(b'_i)$ and $\inter(s)$. Recall that $s$ is included in the path $c_0$, we can then consider the segment
$b'_i$ as a ray which inters $\mathrm{D}_0$ by $y''_i$. By Lemma \ref{SSlmD3}, we know that $b'_i$ must exit $\mathrm{D}_0$ by a point $z''_i$
which is contained in $ s_1\cup\overline{y'_1x'_{j'_1+1}}$, but it would imply that either $\inter(b_i)\cap b_{k+1} \neq \vide$, or
$\inter(b_i)\cap A_1\neq \vide$, which is impossible by assumption. The lemma is then proved. \carre

We can now finish the proof of Lemma \ref{SSlmD2}. Using Lemma \ref{SSlmD1}, we deduce that there exists a sequence of elementary moves which
transforms $\T_1$ into a new triangulation $\T'_1$ containing $A_1$, the edges $b_1,\dots,b_k$, and the segment $\hat{s}$. By replacing
$\overline{x_{j_1}x_{j_1+1}}$ by $\hat{s}$, we get a new maximal tree $A'_1$. Let us show that the number of intersection points of $b_{k+1}$
with $A'_1$, excluding the endpoints of $b_{k+1}$, is at most $m-1$. We have

$$\begin{array}{cl}
  \card\{\inter(b_{k+1})\cap A'_1\}= & \card\{\inter(b_{k+1})\cap A_1\}-\card\{\inter(b_{k+1})\cap\inter(\overline{x_{j_1}x_{j_1+1}})\}+ \\
   & + \card\{\inter(b_{k+1})\cap \inter(\hat{s})\} \\
\end{array}$$

\noindent Let $y$ be a point in  $\inter(b_{k+1})\cap \inter(\hat{s})$, and let $y'=\pi^{-1}_{A_1}(y)$. Let $b'$ be the segment in
$\pi_{A_1}^{-1}(b_{k+1})$ which contains $y'$, note that $y'=b'\cap s$. By Lemma \ref{SSlmD4}, and since $\inter(b')\cap \inter(s_1)=\vide$, it
follows that $b'$ contains a point $z'$ in $\overline{x'_{j'_1}x'_{j'_1+1}}$. We deduce that there is a one-to-one mapping from
$\{\inter(b_{k+1})\cap \inter(\hat{s})\}$ into $\{\inter(b_{k+1})\cap \inter(\overline{x_{j_1}x_{j_1+1}})\}$. Clearly, the point $y_1$ does not
belong to the image of this map, therefore we have

$$\card\{\inter(b_{k+1})\cap\inter(\overline{x_{j_1}x_{j_1+1}})\}\geq \card\{\inter(b_{k+1})\cap \inter(\hat{s})\}+1.$$

\noindent It follows immediately that

$$\card\{\inter(b_{k+1})\cap A'_1\}\leq \card\{\inter(b_{k+1})\cap A_1\}-1=m-1.$$

\noindent The proof of Lemma \ref{SSlmD2} is now complete.  \carre

\noindent From what we have seen, Proposition \ref{SSprD1}, and hence Theorem \ref{SSThD}, follow directly from Lemma \ref{SSlmD2}.\carre

\end{appendices}

\end{document}